%% file: ex_article.tex
\documentclass[hidelinks,onefignum,onetabnum]{siamart220329}

\input{ex_shared}

\begin{document}

\maketitle

% REQUIRED
\begin{abstract}
We establish a theoretical framework of the particle relaxation method for uniform particle generation of Smoothed Particle Hydrodynamics.
We achieve this by reformulating the particle relaxation as an optimization problem.
The objective function is an integral difference between discrete particle-based and smoothed-analytical volume fractions.
The analysis demonstrates that the particle relaxation method in the domain interior is essentially equivalent to employing a gradient descent approach to solve this optimization problem, and we can extend such an equivalence to the bounded domain by introducing a proper boundary term.
Additionally, each periodic particle distribution has a spatially uniform particle volume, denoted as characteristic volume.
The relaxed particle distribution has the largest characteristic volume, and the kernel cut-off radius determines this volume.
This insight enables us to control the relaxed particle distribution by selecting the target kernel cut-off radius for a given kernel function.
\end{abstract}

% REQUIRED

% REQUIRED

\section{Introduction}
Smoothed Particle Hydrodynamics (SPH) \cite{gingold1977smoothed,lucy1977numerical} has gained significant popularity in computational mechanics due to its mesh-free characteristic.
While this characteristic eases the treatments for large deformation and moving boundaries \cite{peng2015sph}, improving the convergence and consistency becomes one of the grand challenges \cite{frontiere2017crksph,vacondio2021grand} because they greatly depend on the particle distribution.
Under irregular particle distributions, the moment constraints are not satisfied \cite{frontiere2017crksph}, which obstructs the zeroth order consistency \cite{litvinov2015towards,basa2009robustness,quinlan2006truncation}.
Moreover, irregular particle distributions compromise the partition of unity and subsequently decrease the convergence rate \cite{litvinov2015towards}.
Consequently, generating uniform particle distributions is a crucial sub-step for SPH, especially for Eulerian SPH with fixed particles \cite{nasar2021high}, but it becomes challenging for domains with complex geometry.
To address this challenge, various strategies have been proposed, such as 
generating particles in lattice structure \cite{dominguez2011development,negi2021algorithms}, 
%or at the center of polygon volume elements, 
employing the weighted Voronoi tessellation method \cite{diehl2015generating}, 
and utilizing the particle relaxation method \cite{litvinov2015towards,fu2019isotropic,zhu2021cad,yu2022level,zhang2023initial}. 
Although these methods have been widely applied in multi-physics simulations \cite{zhang2020sphinxsys,dominguez2022dualsphysics}, they still lack a solid theoretical foundation \cite{negi2021algorithms} because of the difficulties on theoretical treatments of N-body systems.
The particle relaxation method \cite{litvinov2015towards,fu2019isotropic,zhu2021cad,yu2022level,zhang2023initial} offers several advantages, including high efficiency, simplicity, and its body-fitted property, which motivates us to further investigate it theoretically.

The motion of particles in the relaxation method is determined by a pairwise interaction force, which is analogous to swarm dynamics \cite{mogilner2003mutual} in biology. 
In swarm dynamics, 
due to the complexity of studying N-body system with $3N$ degrees of freedom,
a continuous-regime approximation is employed to simplify the analysis, as demonstrated in \cite{bodnar2005derivation,pittayakanchit2016global}. 
A similar strategy was also applied in density functional theory (DFT)  \cite{kohn1996density}
and  dynamical density functional theory (DDFT) \cite{te2020classical} as a tool to investigate crystal structures.
However, these investigations mainly focus on the continuous regime.
Although discrete systems can be represented by the Dirac distribution, a bounded domain requires the introduction of a proper boundary term, and there is no study on how to introduce the boundary term properly.  

We focus on discrete particle systems, and build an optimization framework for particle relaxation, where the objective function is based on  the integral difference between discrete particle-based and smoothed-analytical volume fractions.
The total error is non-negative and diminishes along the individual particle trajectories.
We prove that optimization of the total error via the gradient descent method essentially is equivalent to the particle relaxation method.
Based on this framework, we introduce a boundary term for bounded domains, which is compatible with different relaxed particle distributions.
This optimization framework also provides the capability to predict and control relaxed distributions produced by the relaxation method.
Numerical experiments show full agreement with predictions.
Moreover, beyond the particle generation in SPH, our analysis may also provide a theoretical foundation for other particle-based applications, such as domain decomposition \cite{fu2017physics} and mesh generation \cite{fu2019isotropic,ji2020consistent,ji2021feature}.

The remainder of this article is organized as follows.
We review the original particle relaxation method \cite{zhu2021cad} in Section 2.
Section 3 establishes the optimization framework, and the equivalence between optimizing through the gradient descent method and particle relaxation method is proved.
Section 4 includes the introduction of a boundary term and a comparison with other surface-bounding methods.
In Section 5, we predict relaxed particle patterns based on our optimization framework. 
Finally, we provide conclusions in Section 6. 

\section{Particle relaxation method}
\label{Particle_relaxation_method}
With the particle relaxation method, as detailed in \cite{zhu2021cad}, the domain boundary is firstly defined by a zero level-set function $\phi(\mathbf{x})=0$, followed by a relaxation process.
%In the particle relaxation method \cite{zhu2021cad}, a zero level set $\phi(\mathbf{x})=0$ is generally utilized to describe the domain boundary in advance, followed by a relaxation process. 
The level-set function $\phi(\mathbf{x})$ is a signed distance function with $|\nabla \phi|=1$, and its negative and positive signs respectively indicate the inside and outside domain.
The outward normal vector $\mathbf{n}$ to this boundary is calculated by $\mathbf{n} = \nabla \phi/|\nabla \phi|$.
Furthermore, the level-set function $\phi(\mathbf{x})$ is discretely stored on a Cartesian background mesh with node spacing  $\Delta x$.

During relaxation, with respect to particle $i$, the particle movement is governed by
\begin{equation} \label{original_governing}
	\frac{d\mathbf{u}_i}{dt}=\mathbf{F}_i,
\end{equation}
where $ \mathbf{u}_i $ is the particle velocity, and $ t $ the time. The resulting force $\mathbf{F}_i$ is defined as the summation of pairwise repulsive forces
%\math
\begin{equation} \label{repulsive_force}
	\mathbf{F}_i = -\frac{2}{m_i}\sum_j{p_0V_iV_j\nabla_i W(\Vert\mathbf{x}_i-\mathbf{x}_j\Vert,h)},
\end{equation}
where $ m_i $ is the particle mass, $ p_0 $ the pressure, $ V $ the particle volume, and $ h $ the cut-off radius. The subscript $ j $ denotes the neighbor particle. Here, $W$ represents the $5^{th}$-order Wendland function, and $p_0$ is typically set to $1$. 
Additionally, due to the common utilization of $m_i = 1$ and $V_i=V_j=(\Delta x)^d$, where $d$ denotes the dimension, $ \mathbf{F}_i $ becomes equivalent to the particle acceleration.

Hence, the particle position is updated by 
\begin{equation}\label{eq:updatezhu}
	\mathbf{r}^{n+1}_i=\mathbf{r}^{n}_i+\frac{1}{2}\mathbf{F}_i^{n}\Delta t^2,
\end{equation}
%where $\Delta t \leq 0.25 \sqrt\frac{h}{\max|\mathbf{F}_i^{n}|}$.
where $\Delta t \leq 0.25 \sqrt{h/\max|\mathbf{F}_i^{n}|}$.

However, Eq. (\ref{eq:updatezhu}) only considers the displacement resulting from particle acceleration.
This implies the necessity of resetting the initial velocity $\mathbf{u}_i^{n}$ to zero at each time step \cite{zhu2021cad}. 
By adopting the concept of pseudo force \cite{bernoff2013nonlocal}, we can equivalently replace the displacement $ \frac{1}{2}\mathbf{F}_i^{n}\Delta t^2$ with $ \mathbf{U}_i^{n}\Delta T $, eliminating the need for the previously mentioned velocity resetting.

Consequently, in present work, the particle position could also be updated by
\begin{equation}
	\mathbf{r}^{n+1}_i=\mathbf{r}^{n}_i + \mathbf{U}_i^{n} \Delta T,
	\label{eq:stokes}
\end{equation}
where the velocity $\mathbf{U}_i^{n}$ is proportional to the corresponding pseudo force exerted on particle $ i $.
%$\Delta T = \frac{ch}{{\max}|\mathbf{U}_i^{n}|}$ is the time step,
$\Delta T = {ch/{\max}|\mathbf{U}_i^{n}|}$ is acoustic time step,
and the coefficient $c$ can control the position increment $ \Vert \Delta \mathbf{r} \Vert_{2,\infty} $ through
\begin{equation}
	\Vert \Delta \mathbf{r} \Vert_{2,\infty} = \Vert \mathbf{r}^{\,n+1} - \mathbf{r}^{\,n}\Vert_{2,\infty} = \max|\mathbf{U}_i^{n}| %\frac{ch}{{\max}|\mathbf{U}_i^{n}|}
	\Delta T=ch.
	\label{eq:controlsize} 
\end{equation}

Note that, corresponding to different $ \Delta t $ in Eq. (\ref{eq:updatezhu}), there exists a unique $c$ to maintain the equivalence between Eq. (\ref{eq:updatezhu}) and (\ref{eq:stokes}). For example, $c = 1/32$ corresponds to  $\Delta t = 0.25\sqrt{h/\max|\mathbf{F}_i^{n}|}$.

\section{Equivalence of optimization and particle relaxation}
\label{chap:formula}

In this section, we prove that in the domain interior, the particle relaxation method as detailed in \cite{zhu2021cad} is equivalent to optimization of a total error by the gradient descent approach. 
The total error is defined by integrating an error density, namely, the integral difference between discretized particle-based and smoothed-analytical volume fractions.

The interior domain can be treated as a special case of bounded domain by a negative signed distance, i.e. $\forall \mathbf{x} \in \mathbb{R}^d, \ \phi(\mathbf{x})\equiv -\infty$.
Consequently, we build a framework for bounded domains.
We assume that domain $\mathscr{D}=\{\mathbf{x} \in \mathbb{R}^d,\ \phi(\mathbf{x}) \leq 0\}$ has a smooth boundary $\partial \mathscr{D}$,
total volume $\mathcal{V}$,
and $N$ particles are located at $\{\mathbf{x}_i\}_{i=1}^N \subset \mathscr{D}$.
Moreover, for a given kernel cut-off radius $h$, we suppose that there exists an upper bound  $M$ for the number of neighbor particles located within the kernel support surrounding each particle.
Also, the number of near-boundary particles in $\mathscr{D}_b = \{\mathbf{x} \in \mathbb{R}^d,\  \phi(\mathbf{x}) \in (-h,0)\}$ and inner particles in $\mathscr{D}_i = \{\mathbf{x} \in \mathbb{R}^d,\  \phi(\mathbf{x}) <h\}$ respectively are on the order of $O(N^{(d-1)/d})$ and $O(N)$.

\subsection{Expression for particle-based and smoothed-analytical volume fractions}
According to the normalization property of the kernel function $W$, when considering particles at different position $\mathbf x \in \mathscr{D}$, the integration of $W$ within the support domain $\Omega$ can be regarded as a volume fraction at $\mathbf{x}$
\begin{equation}\label{normalization}
	0\leq\int_{\Omega}{W(\Vert\mathbf{x}-\mathbf{x'}\Vert,h)d\mathbf{x'}}\leq 1.
	%=
	%\left \{
	%\begin{split}
	%	&1, &\Omega \subset \mathscr{D} \\
	%	&(0,1), &\Omega \cap \partial \mathscr{D}\neq\emptyset \\
	%	&0, &\Omega \cap \mathscr{D}=\emptyset
	%\end{split}
	%\right ..
\end{equation}
Here the left and right equalities hold if and only if 
$\Omega \cap \mathscr{D}=\emptyset$ and $\Omega \subset \mathscr{D}$, respectively.
Upon discretization with particles, the particle-based volume fraction can be written as 
\begin{equation}
	\alpha(\mathbf{x}) = \sum_{j}{W(\Vert \mathbf{x}-\mathbf{x}_j \Vert,h)v_0},
	\label{vf}
\end{equation}
where the real average volume $v_0= \mathcal{V}/N$ can maintain the consistency of total volume $\int_{\mathbb{R}^d}\alpha(\mathbf{x})d\mathbf{x}=\mathcal{V}$ exactly.
Because the particles are initialized at the cell center with a small random shift, we have
\begin{equation}
\lim_{N \to \infty} v_0/(\Delta x)^d = \lim_{N \to \infty} \mathcal{V}/[N(\Delta x)^d]=1.
\end{equation} 
Therefore for sufficient large $N$, we can still approximate the average particle volume $v_0$ as the cell volume $(\Delta x)^d$.

The characteristic function of a domain is a sharp volume fraction, and it can be written as a composite function of the Heaviside function and level-set function \cite{osher2005level}
\begin{equation}
	H(-\phi(\mathbf{x}))=\left\{
	\begin{aligned}
		1,\quad \mathbf{x} \in \Omega\\
		0,\quad \mathbf{x}\not\in \Omega
	\end{aligned}
	\right..
	\label{eq:charac}
\end{equation}
Thus the smoothed-analytical volume fraction can be obtained by smoothing it
\begin{equation}\label{smoothcharac}
	0\leq
	P(-\phi(\mathbf{x}))\\
	=1-2\int_{-\phi(\mathbf{x})}^{h}{W_1(\max\{r,0\},h)dr}
	\leq 1,
\end{equation}
where $W_1(r,h) = AW(r,h)$, and A is constant satisfying $\int_{-h}^{h}AW(r,h)dr=1$.
The left and right equality signs in the Eq. (\ref{smoothcharac}) attain if and only if $\phi(\mathbf{x}) \in [0,\infty)$ and $\phi(\mathbf{x})\in (-\infty,-h]$, respectively.

While the particle-based volume fraction Eq. ({\ref{vf}}) depends on the particle distribution,
the smoothed-analytical volume fraction Eq. (\ref{smoothcharac}), remains independent of the particle distribution.
The smoothed-analytical volume fraction can be regarded as an ideal volume fraction, and the difference between these two kinds of volume fraction serves as a measure of the quality of particle distribution.
\subsection{Definition of error density and total error}
The error density $e(\mathbf{x})$ is defined as the difference between the particle-based and smoothed-analytical volume fractions Eq. (\ref{vf}) and Eq. (\ref{smoothcharac}),
\begin{equation}
	\begin{split}
		e(\mathbf{x}) 
		&= \alpha(\mathbf{x})-P(-\phi(\mathbf{x}))\\
		&=\sum_{j}{W(\Vert \mathbf{x}-\mathbf{x}_j \Vert,h)v_0} -1+ 2\int_{-\phi(\mathbf{x})}^{h}{W_1(\max\{r,0\},h)dr},
	\end{split}
\end{equation}
and the total error $E$ can be measured by integrating the error density $e(\mathbf{x})$
\begin{align}
	E&=\int_{\mathscr{D}}{e(\mathbf{x})d\mathbf{x}}\notag\\
	&= \int_{\mathscr{D}}\sum_{j}{W(\Vert \mathbf{x}-\mathbf{x}_j \Vert,h)v_0}d\mathbf{x}-\mathcal{V}+2\int_{\mathscr{D}}\int_{-\phi(\mathbf{x})}^h{W_1(\max\{r,0\},h)dr}d\mathbf{x}\label{eq:conbcf}\\
	&\approx \sum_{i,j}{W(\Vert \mathbf{x}_i-\mathbf{x}_j \Vert,h)v_0 \tilde{v}}-\mathcal{V}+\sum_{i}\left[2\int_{-\phi(\mathbf{x}_i)}^h{W_1(\max\{r,0\},h)dr}\right]\tilde{v}\label{eq:disbcf}\\
	&=\sum_{i,j}{W(\Vert \mathbf{x}_i-\mathbf{x}_j \Vert,h)v_0 \tilde{v}}-\mathcal{V}+O(h),\label{eq:disesti}
\end{align}
Here, the supremum of the average particle volume $\tilde{v} = \sup_{(\mathbf{x}_1,\dots,\mathbf{x}_n)} \{\bar{v}(\mathbf{x}_1,\dots,\mathbf{x}_n)\}$ is used to discretize the domain $\mathscr{D}$ in order to ensure the positivity of total error, as detailed in Section \ref{chap:positive}.
The average particle volume is defined as $\bar{v}=\sum_p{v_p}/N$.

The last equality holds because the last term in Eq. (\ref{eq:disbcf}) can be estimated as $O(h)$ as follows.
The volume $v_i$ of each particle $i$ is bounded by
\begin{equation}
	\frac{1}{W(0,h)}\geq v_i = \frac{1}{\sum_j W(\Vert \mathbf{x}_i- \mathbf{x}_j\Vert,h)} \geq \frac{1}{MW(0,h)}.\label{eq:volcontrol}
\end{equation}
Here, the left inequality in Eq. (\ref{eq:volcontrol}) results from the positivity of the kernel function and $W(\Vert \mathbf{x}_i- \mathbf{x}_i\Vert,h)=W(0,h)$.
The right inequality results from $W(0,h) \geq W(r,h)$ for arbitrary $r$.
According to the definition of kernel function, both sides in Eq. (\ref{eq:volcontrol}) have  order of $h^d$.
Thus, we have the following estimate
\begin{equation}
	v_i 
	\propto \bar{v} 
	\propto \tilde{v} 
	\propto v_0 \propto h^d \propto \frac{1}{N}.\label{volumeestimation}
\end{equation}
By multiplying $\tilde{v}=O(h^d)$ with the number of boundary particles and the bounded integral $ 0 \leq \int_{-\phi(\mathbf{x}_i)}^h{W_1(\max\{r,0\},h)dr} \leq 1/2$, the order of the last term in Eq. (\ref{eq:disbcf}) is estimated as $O(h)$.

%Note that the last term in Eq. (\ref{eq:conbcf}) also has order $O(h)$ because the inner integral is is bounded by $[0,1]$ and only non-zero in $\mathscr{D}_b$, whose volume is $O(h)$. 

\subsection{Non-negative property of the total error}
\label{chap:positive}
The non-negative property of the total error $ E $ can be proved as follows
\begin{align}
	E&=\sum_i\sum_j{W(\Vert \mathbf{x}_i -\mathbf{x}_j \Vert,h)v_0\tilde{v}}-\mathcal{V}
	+O(h)\label{eq:dissipative_original}\\
	&\geq \left(\sum_i\sum_j{W(\Vert \mathbf{x}_i -\mathbf{x}_j \Vert,h)}v_0\tilde{v}\right)-\mathcal{V}\label{eq:removeOh}\\
	&= \left(\sum_i{\frac{1}{v_i}}v_0 \tilde{v}\right)-\mathcal{V}\\
	&= \left( \frac{\mathcal{V}\tilde{v}}{H}\right)-\mathcal{V}\\
	&\geq \mathcal{V}\left(\frac{\tilde{v}}{\bar{v}}-1\right) \geq 0, \label{meanineq}
\end{align}
where $H=N/(\sum_i{1/v_i})$ is the harmonic mean.

Note that, Eq. (\ref{eq:removeOh}) is obtained by neglecting the non-negative term $O(h)$ in Eq. (\ref{eq:dissipative_original}), and Eq. (\ref{meanineq}) is derived from the mean inequality $\bar{v} \geq H$.
Furthermore, as $E$ is always non-negative, $E = |E|$, and the absolute error can be minimized by optimizing the total error.

\subsection{Equivalence to particle relaxation in domain interior}
\label{chap:equi}
The last term in Eq. (\ref{eq:disbcf}), i.e. the boundary term vanishes in domain interior.
When optimizing the total error via the gradient descent method, the negative gradient of the total error with respect to particle $i$ is written as
\begin{equation}
	-\frac{\partial E}{\partial \mathbf{x}_i} \approx
	%2\sum_{j \neq i} W'(\Vert \mathbf{x}_i-\mathbf{x}_j\Vert,h)\frac{\mathbf{x}_i-\mathbf{x}_j}{\Vert \mathbf{x}_i-\mathbf{x}_j\Vert}v_0\tilde{v}=
	2\sum_{j \neq i} W'(\Vert \mathbf{x}_i-\mathbf{x}_j\Vert,h)\mathbf{e}_{ij}v_0\tilde{v},
	\label{eq:ngws}
\end{equation} 
where $\mathbf{e}_{ij} = (\mathbf{x}_{j}-\mathbf{x}_{i})/\Vert\mathbf{x}_{j}-\mathbf{x}_{i}\Vert$, and $j \neq i$ sums over all neighbors of particle $i$.

As $-\partial E/\partial \mathbf{x}_i$ is parallel to $\mathbf{F}_i^{n}$ in Eq. (\ref{repulsive_force}) and $\mathbf{U}_i^{n}$ in Eq. (\ref{eq:stokes}),
minimizing the total error along the negative gradient is equivalent to particle relaxation \cite{zhu2021cad}.
We stress that the total error diminishes along the individual particle trajectories, as 
\begin{equation}
	\frac{dE}{dt} = \sum_{i=1}^N \frac{\partial E}{\partial \mathbf{x}_i} \cdot \mathbf{U}_i^{n}\propto-\nabla E \cdot \nabla E \leq 0.
\end{equation}

For bounded domains, the last term in Eq. (\ref{eq:disbcf}) does not vanish and thus the equivalence does not hold strictly.
%However, this term serves a similar purpose to the surface bounding method \cite{zhu2021cad} and has more attracting properties, which is detailed in Section \ref{chap:boundedregion}.
%
\section{Extension to bounded domain}
\label{chap:boundedregion}
\subsection{The boundary term}
For a bounded domain, because of incomplete support the resulting repulsive force, i.e. $\sum_jW'(\Vert \mathbf{x}_i-\mathbf{x}_j\Vert,h)\mathbf{e}_{ij}v_0\tilde{v}$ as mentioned in \cite{zhu2021cad}, propels boundary particles across the boundary.
Consequently, a compensating force is essential to counteract this effect and prevent the boundary particles from escaping the domain.
\begin{figure}[H]
	\centering
	\includegraphics[width = 0.8\textwidth]{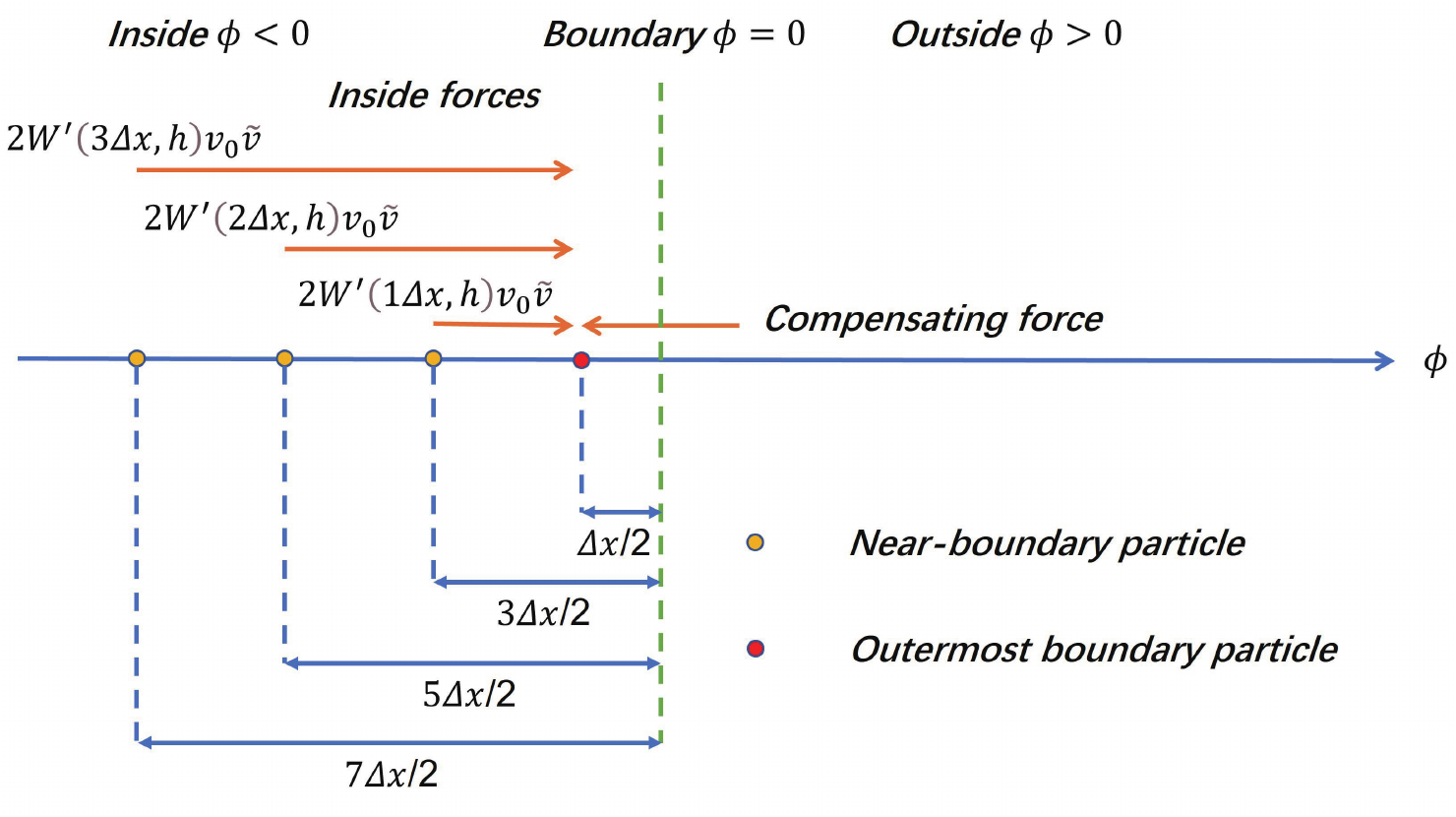}
	\caption{
		%The force introduced by Eq. (\ref{bcf}) should balance the inner force at least in 1D. 
		The definition of the compensating force for boundary particles in 1D. The uniform particle spacing is $\Delta x$.
	}
	\label{fig:balance}
\end{figure}

Generalizing the case in Fig. \ref{fig:balance}, the magnitude of the resulting force exerted on $l^{th}$-layer boundary particle at $\phi=-(2l-1)\Delta x/2$ is
\begin{equation}
	f\left(\frac{2l-1}{2}\Delta x\right)=-2\sum_{k=l}^{\infty}{W_1'(k\Delta x,h)v_0\tilde{v}}.
\end{equation}
Here, the resulting force points towards the outward normal direction and can be further approximated as
\begin{equation}
	-2\sum_{k=l}^{\infty}{W_1'(k\Delta x,h)v_0\tilde{v}} \approx -2\tilde{v}\int_{\frac{2l-1}{2}\Delta x}^{h}{W_1'(r,h)dr}=2W_1\left(\frac{2l-1}{2}\Delta x,h\right)\tilde{v}.
	\label{integral}
\end{equation}

For particles outside the domain, i.e., $\phi(\mathbf{x})>0$, a compensating force is necessary, which should be as large as possible but remain continuous with $2W_1(-\phi(\mathbf{x}),h)\tilde{v}$ at $\phi(\mathbf{x}) = 0$.
Consequently, Eq. (\ref{integral}) is modified as to be applicable to both particles inside and outside the domain
\begin{equation}
	f(-\phi(\mathbf{x})) = 2W_1(\max\{-\phi(\mathbf{x}),0\},h)\tilde{v}.
	\label{eq:generalbf}
\end{equation}

Compared to Eq. (\ref{eq:ngws}) for the domain interior, as mentioned in Section \ref{chap:equi}, the negative gradient of total error $ E $ for the bounded domain is
\begin{equation}
	-\frac{\partial E}{\partial \mathbf{x}_i} \approx
	2\sum_{j \neq i} W'(\Vert \mathbf{x}_i-\mathbf{x}_j\Vert,h)\mathbf{e}_{ij}v_0\tilde{v}-2W_1(\max\{-\phi(\mathbf{x}_i),0\},h)\tilde{v}\mathbf{n}.
	\label{eq:ngb}
\end{equation}
Here, the magnitude of the last term is same as that of Eq. (\ref{eq:generalbf}), which implies that our formulation can balance the repulsive forces near the boundary. This term can be regarded as a boundary correction for particle relaxation in bounded domains.

As we use an analytical form to replace the summation in Eq. (\ref{integral}) for a bounded domain, ensuring the accuracy of the approximation requires a relatively large kernel cut-off radius $h>2.0\Delta x$. The supporting numerical tests are provided in Appendix \ref{validationpatt}.
\subsection{Comparison with the surface bounding method}
\begin{figure}[H]
	\centering
	\includegraphics[width = 0.75\textwidth]{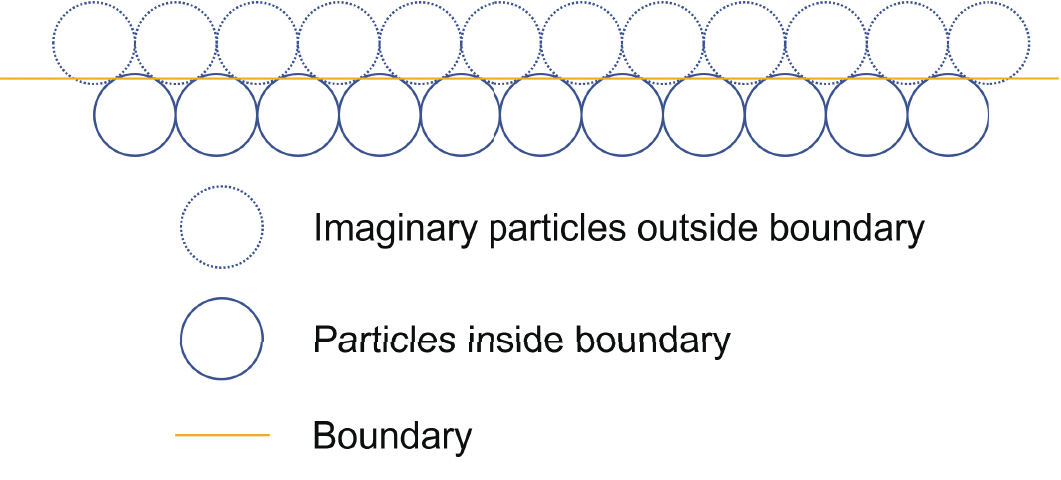}
	\caption{Hexagonal distribution, a counter example for \cite{zhu2021cad,zhang2023initial}.}
	\label{fig:straightboundary}
\end{figure}

\begin{figure}[H]
	\centering
	\includegraphics[width = 0.75\textwidth]{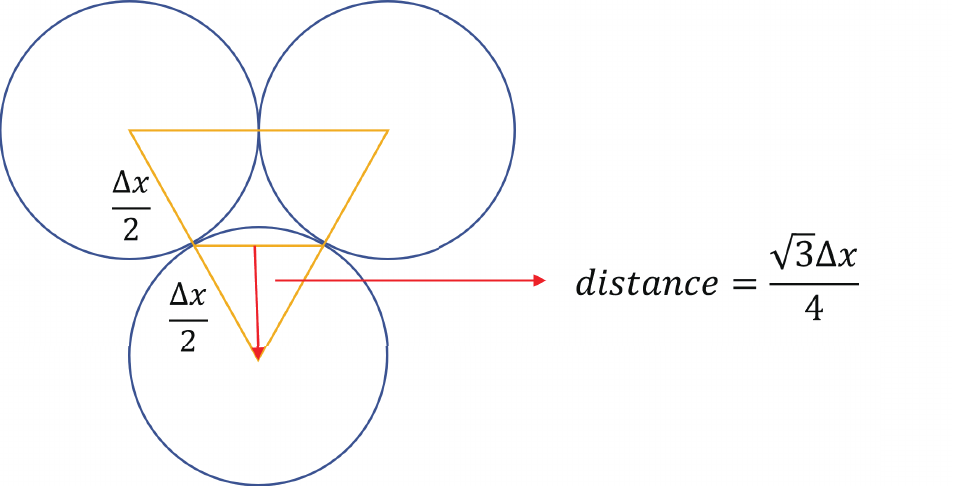}
	\caption{Distance calculation of the first-layer particles to the boundary.}
	\label{fig:distance}
\end{figure}

\begin{figure}[H]
	\centering
	\includegraphics[width = 0.75\textwidth]{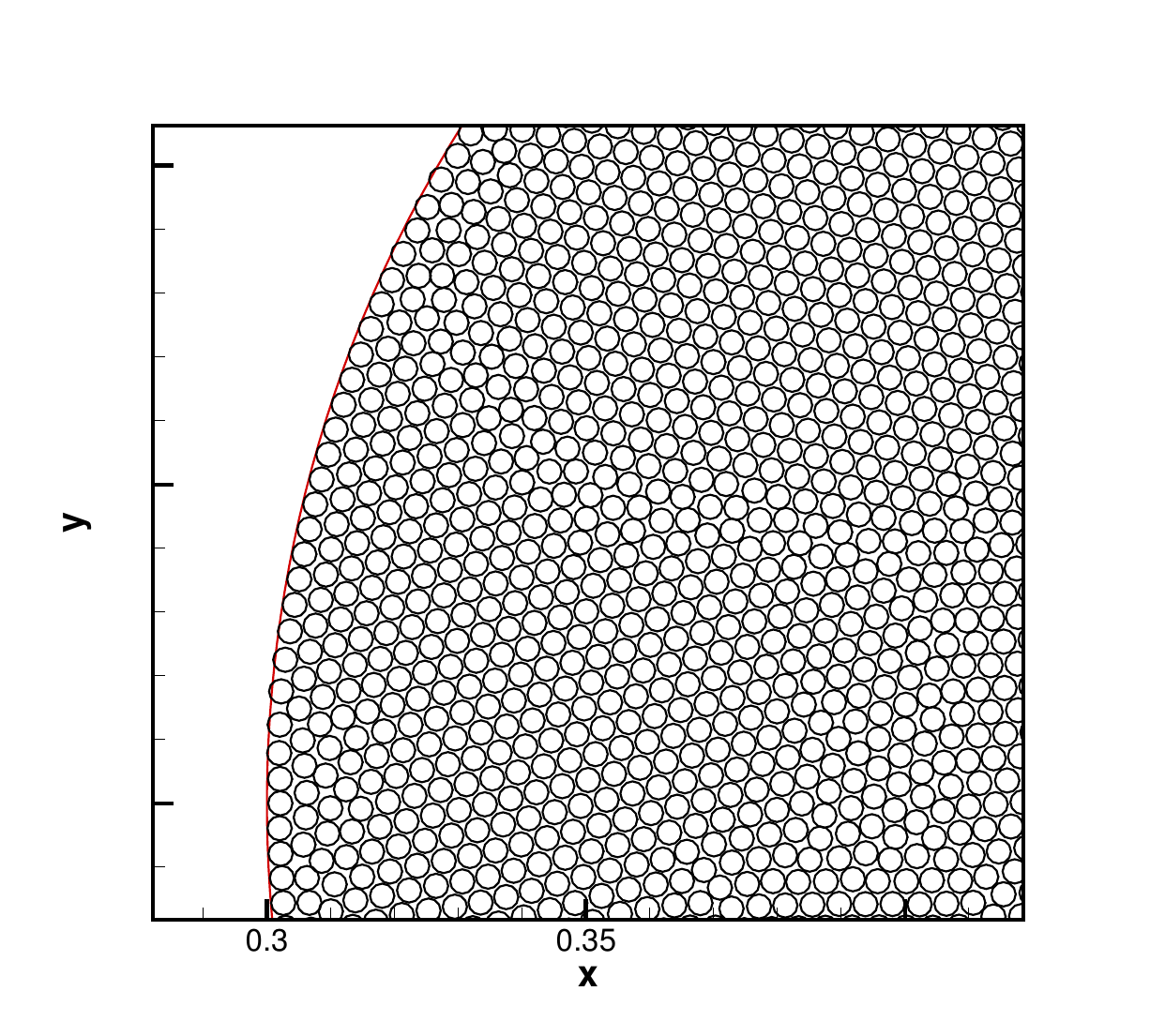}
	\caption{A hexagonal distribution inside a circle obtained by the present method.}
	\label{fig:countereg}
\end{figure}

The surface bounding method in \cite{zhu2021cad} is defined as
\begin{equation}
	\vec{r}_i \leftarrow \left\{
	\begin{aligned}
		&\vec{r}_i-\left(\phi(\mathbf{x})+\frac{1}{2}\Delta x\right)\vec{n}, &\phi(\mathbf{x})>-\frac{1}{2}\Delta x\\
		&\vec{r}_i,& else
	\end{aligned}
	\right.,
\end{equation}
where '$\leftarrow$' implies assigning the value on the right-hand side to the left-hand side.
This method only corrects outermost particles with $\phi(\mathbf{x})>-\Delta x/2$, which leads to a non-uniform distribution among several near-boundary layers.
This method was improved in \cite{zhang2023initial} by identifying different particle layers and aligning particles to corresponding layer at $\phi=-(2l-1)\Delta x/2$, where $l$ is the index of the layer.
However, imposing the particles at a specific distance from  boundary, i.e. an odd integer times of $\Delta x/2$, is not entirely appropriate because it applies only to square/cubic lattices.
A counterexample is the case of a hexagonal distribution adjacent to a straight line boundary, Fig. \ref{fig:straightboundary}, where the first-layer particles should be positioned at a distance of $\sqrt{3}\Delta x/4$ from the boundary, as calculated 
in Fig. \ref{fig:distance}.
Similarly, setting the distance from the $l^{th}$ layer as $(2l-1)\Delta x/2$ also is inappropriate.
In contrast to \cite{zhu2021cad,zhang2023initial}, our formulation introduces a boundary correction Eq. (\ref{eq:generalbf}), which naturally confines the particle to the given domain.
This term provides a suitable confining force at any position within $\mathscr{D}_b$, making our method compatible to various relaxed particle distributions, as illustrated in Fig. \ref{fig:countereg}.
A straightforward application is to generate particles in two adjacent regions independently, as detailed in Appendix \ref{chap:nest}.
Moreover, unlike the static confinement method \cite{yu2022level} our method is parameter-free and does not rely on other variables, such as the volume fraction of cut cells.
\section{Relaxed particle distribution as target patterns}
\subsection{Candidate particle distributions for domain interior}
\begin{figure}[tb!]
	\centering
	\includegraphics[width = 1.0\textwidth]{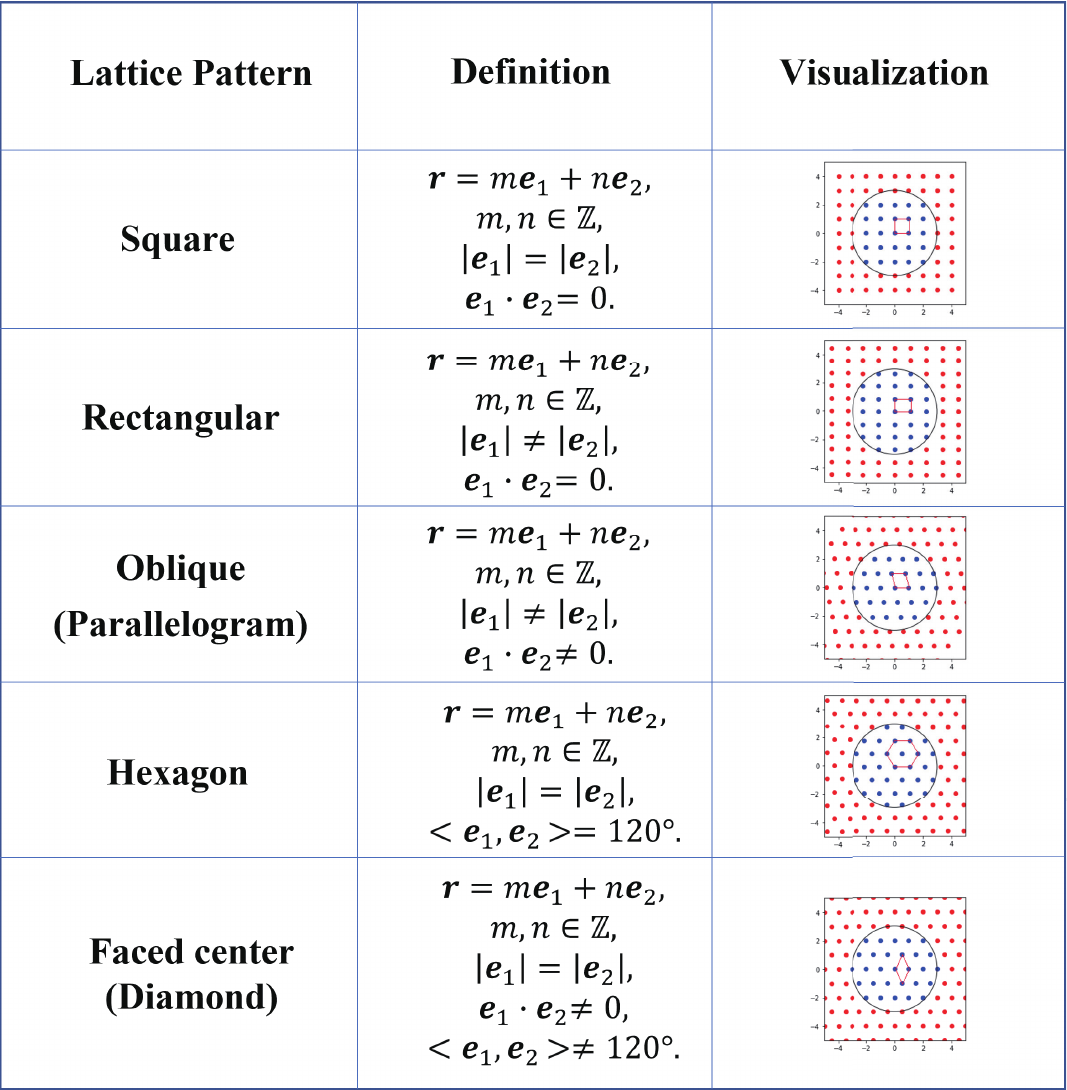}
	\caption{List of all 2D Bravais Lattices
		, see also \cite{kittel2005introduction}.}
	\label{tab:2D-lattices}
\end{figure}

In 2D unbounded domains, there are precisely five types of periodic lattice distribution patterns, known as Bravais lattices \cite{kittel2005introduction}, as listed in Fig. \ref{tab:2D-lattices}.
They respectively correspond different total errors $E$.
Based on the definition of $E$ in Eq. (\ref{eq:disbcf}), only the first term is influenced by the particle distribution in the domain interior
\begin{equation}
	E_P  = \left(\sum_i{\frac{1}{v_i}}v_0 \tilde{v}\right).
	\label{eq:ep}
\end{equation}

For periodic distributions, the volume of all particles is constant, i.e., $v_i \equiv v_c$, denoted by characteristic volume. $E_P$ can be further rewritten as
\begin{equation}
	E_P  = 
	\mathcal{V}\frac{\tilde{v}}{v_c}.
	\label{eq:periodicwhole}
\end{equation}
Here, for a specified kernel function, $v_c$ is determined by kernel-support cut-off radius $h$, thus $E_P $ essentially varies with $h$.
Hence, the predicted pattern of particle distribution, corresponding to the smallest $E_P $,  also changes with varying $h$.  
Tab. \ref{tab:sum_prediction} provides the predictions for the range of $h/\Delta x$ from 1.5 to 3.5.

%\begin{table}[tb!]
%	\centering
%	\caption{2D Bravais Lattices}
%	\label{tab:2D-lattices}
%	\begin{tabular}{|m{5cm}<{\centering}|m{5cm}<{\centering}|m{5cm}<{\centering}|}
	%		\hline
	%		\textbf{Lattice Name} & \textbf{Lattice Points} & \textbf{Lattice Visualization} \\ \hline
	%		Square & $\mathbf{r}=m\mathbf{e}_1+n\mathbf{e}_2$, $m,n \in \mathbb{Z}$, where $|\mathbf{e}_1|=|\mathbf{e}_2|$, and $\mathbf{e}_1 \cdot \mathbf{e}_2=0$ always hold.& \includegraphics[scale=0.35]{Fig/SquareLattice.png} \\ \hline
	%		Rectangular & $\mathbf{r}=m\mathbf{e}_1+n\mathbf{e}_2$, $m,n \in \mathbb{Z}$, where $|\mathbf{e}_1|\neq|\mathbf{e}_2|$, and $\mathbf{e}_1 \cdot \mathbf{e}_2=0$ always hold.& \includegraphics[scale=0.35]{Fig/RectangleLattice.png} \\ \hline
	%		Oblique(Parallelogram) & $\mathbf{r}=m\mathbf{e}_1+n\mathbf{e}_2$, $m,n \in \mathbb{Z}$, where $|\mathbf{e}_1|\neq|\mathbf{e}_2|$, and $\mathbf{e}_1 \cdot \mathbf{e}_2 \neq 0$ always hold.& \includegraphics[scale=0.35]{Fig/ParallelogramLattice.png} \\ \hline
	%		Hexagonal & $\mathbf{r}=m\mathbf{e}_1+n\mathbf{e}_2$, $m,n \in \mathbb{Z}$, where $|\mathbf{e}_1|=|\mathbf{e}_2|$, and $<\mathbf{e}_1,\mathbf{e}_2>=120^{\circ}$ always hold.& \includegraphics[scale=0.35]{Fig/HexagonLattice.png} \\ \hline
	%		Centered rectangular(Diamond)
	%		 & $\mathbf{r}=m\mathbf{e}_1+n\mathbf{e}_2$, $m,n \in \mathbb{Z}$, where $|\mathbf{e}_1|\neq|\mathbf{e}_2|$, and $\cos<\mathbf{e}_1,\mathbf{e}_2>=|\mathbf{e}_1|/|2\mathbf{e}_2|$ always hold.& \includegraphics[scale=0.35]{Fig/DiamondLattice.png} \\ \hline
	%	\end{tabular}
%\end{table}

\begin{table}[tb!]
	\centering
	\caption{Predicted 2D lattice pattern as function of kernel cut-off radius}
	\begin{tabular}{|c|c|}
		\hline
		$h/\Delta x$ &Pattern\\ \hline
		$[1.5,1.8]$, $[2.25,2.75]$, $[3.3,3.5]$ & Hexagonal distribution \\ \hline
		$[2.1,2.2]$, $[2.8,3.25]$ & Parallelogram distribution \\ \hline
		$[1.85,2.05]$ & Square distribution \\ \hline
	\end{tabular}
	\label{tab:sum_prediction}
\end{table}

\subsection{Extension to bounded domains}
The term $E_P$ of bounded domains additionally includes the last term in Eq. (\ref{eq:disbcf})
\begin{equation}
	E_P = \left(\sum_i{\frac{1}{v_i}}v_0 \tilde{v}\right) + \sum_{i}\left[2\int_{-\phi(\mathbf{x}_i)}^h{W_1(\max\{r,0\},h)dr}\right]\tilde{v}.
\end{equation}
However, the second term is $O(h)$, and thus becomes negligible as $N$ increases.
Consequently, the present $E_P$ can be simplified as Eq. (\ref{eq:ep}).
Although near-boundary particles in bounded domains have larger volume,
conclusions from domain interior still hold in bounded domains, which can be proved as follows.

We consider a periodic distribution $d \in \mathcal{D}$ with characteristic volume $v_c^d \propto 1/N$,  $\mathcal{D}$ represents the set of all possible candidate distributions.
It is assumed that there are $N_1=O\left(N^{(d-1)/d}\right)$ near-boundary particles, each with unknown volume but bounded by $v_b=O(h)$.
Moreover $N_2=O(N)$ particles in the domain interior with characteristic volume $v_c^d$. 
Given that boundary particles possess larger volumes than those in the interior due to incomplete support, $E_P$ can be further divided into two parts
\begin{align}
	E_P 
	&= \sum_i \sum_j W_{ij}\tilde{v}v_0\notag\\
	&= \sum_i \frac{1}{v_i}\tilde{v}v_0\notag\\
	&= \sum_p \frac{1}{v_p}\tilde{v}v_0 +\sum_q \frac{1}{v_q}\tilde{v}v_0,
\end{align}
where $p$ sums over the $N_1$ near-boundary particles and $q$ sums over the $N_2$ inner particles.
Because the volume of a particle $p$ is bounded by $v_b$, we can find the lower and upper bounds of $E$ as
\begin{equation}
	\tilde{v}v_0\left(-N_1 \frac{1}{v_b} + N_2 \frac{1}{v_c^d}\right) \leq E_P \leq  \tilde{v}v_0\left(N_1 \frac{1}{v_b} + N_2 \frac{1}{v_c^d}\right).
\end{equation}
Apparently, the two sides have the same limit $N_2\tilde{v}v_0/v_c^d$ as $N \to \infty$, because 
\begin{equation}
\tilde{v}v_0 N_1/v_b =O(N^{-1/d}),\quad \tilde{v}v_0 N_2/v_c^d = O(1).
\end{equation}
Thus, $E_P$ is determined by the characteristic volume $v_c^d$ of the pattern $d$ in the domain interior.

Numerically, we employ a 2D case to verify our analysis.
We consider a circular domain centered at $(0.5,0.5)$ with a radius of $r=0.2$, and a total volume of $\mathcal{V} = \pi r^2$.
There is a $256 \times 256$ background mesh within the computational domain $[0,1] \times [0,1]$ and 8186 particles in the circle.
We define an effective distance of particles as $dp=\sqrt{\mathcal{V}/N}=\sqrt{v_0}$
%and the corresponding effective cut-off radius of kernel function as
%$h_p = (h/dp)\Delta x$.
We find that the relation between $(h/dp)$ and the corresponding patterns obeys Tab. \ref{tab:sum_prediction}.
With this relation, we can control the relaxed particle distribution as a target pattern by selecting a suitable kernel cut-off radius.
This insight may also be applied in other applications such as \cite{fu2019isotropic,ji2020consistent,ji2021feature}.
The details of the simulation results are shown in the Appendix \ref{implementation}.

In addition, we conducted tests in a 3D cubic domain with periodic boundaries, finding that the optimal isotropic pattern, specifically a cubic lattice, is only achieved when $h=2\Delta x$.
This observation suggests that for the $5^{th}$-order Wendland kernel function, it is not feasible to have a unique cut-off radius $h \in [2.0\Delta x, 3.5\Delta x]$, which achieves best isotropic property in both 2D and 3D domains.

\section{Comparison with other particle relaxation method}
Litvinov et.al \cite{litvinov2015towards} introduced an inertia term and a viscous term with viscosity $\mu$, which can be summarized as
\begin{equation}
	\left\{
	\begin{aligned}
		\mathbf{a}_i^{n+1}&=2\sum_{j}W'(\Vert \mathbf{x}_i -\mathbf{x}_j \Vert,h)\mathbf{e}_{ij}v_0^2
		+\mu \frac{\mathbf{u}_i-\mathbf{u}_j}{\Vert \mathbf{x}_i -\mathbf{x}_j \Vert}W'(\Vert \mathbf{x}_i -\mathbf{x}_j \Vert,h)v_0^2\\
		\mathbf{r}_i^{n+1}&=\mathbf{r}_i^{n}+\mathbf{u}_i^n \Delta t+\frac{1}{2}\mathbf{a}_i^{n}\Delta t^2\\
		\mathbf{u}_i^{n+1}&= \mathbf{u}_i^{n}+\mathbf{a}_i^{n}\Delta t
	\end{aligned}
	\right..
	\label{sergeyrelax}
\end{equation}
%where $V= \mathcal{V}/N$ is the real average particle volume.

\begin{figure}[H]
	\centering
	\includegraphics[width = 0.6\textwidth]{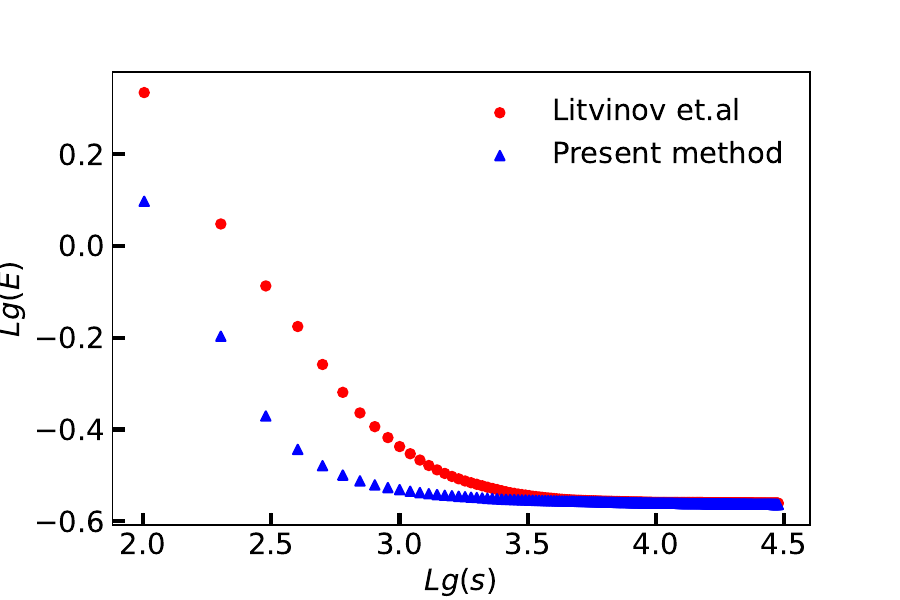}
	\caption{The decay rate of present relaxation method \cite{zhu2021cad} and the one mentioned in \cite{litvinov2015towards}, $s$ is the number of steps, and $h=2.6\Delta x$.
	}
	\label{fig:energycomp}
\end{figure}

Neither the present method i.e. Eq. (\ref{eq:stokes}) nor the method Eq. (\ref{sergeyrelax}) gives periodic particle distributions at early and intermediate stages , as reported in \cite{litvinov2015towards}.
However, we find that with sufficient number of iterations, our method and that of \cite{litvinov2015towards} yield identical periodic particle distribution at late stages.
Moreover, as shown in Fig. \ref{fig:energycomp}, the present particle relaxation \cite{zhu2021cad} demonstrates higher efficiency at the early stages.
This is attributed to fact that the relaxation method always aligns the negative gradient of the total error, in contrast to \cite{litvinov2015towards}, which deviates from the negative gradient of the total error.
\section{Conclusion}
We propose an optimization framework for systematic particle relaxation.
Optimization objective is the total error, defined as the integral difference between a discrete particle-based volume fraction and smoothed-analytical volume fraction.
Minimization of the total error by gradient descent is equivalent to the particle relaxation method in the domain interior.
Near boundaries, we introduce a proper boundary correction, which prevents particles from escaping the specified domain.
Compared to the surface bounding methods \cite{zhu2021cad,zhang2023initial}, this correction also results in a better compatibility between boundary and interior particle distribution.
The effectiveness of the method has been demonstrated by numerical results.
We have shown that for a given kernel function the relaxed particle distribution is determined by the kernel cut-off radius.
Thus, targeted particle distributions are controlled by selecting the cut-off radius and the kernel function.

Beyond generating uniform particle generation for SPH, the present analysis may be applied for other applications, such as mesh generation \cite{fu2019isotropic,ji2020consistent,ji2021feature} and domain decomposition \cite{fu2017novel}.
There are also some further challenges that needed to be solved.
These include applying particle relaxation methods confined to co-dimension one manifolds and constructing kernel function with a unified cut-off radius that optimizes isotropic properties in  both 2D and 3D domains.

\appendix
\section*{Appendix}
\section{Implementation}
\label{implementation}
In our implementation, following the definition Eq. (\ref{eq:stokes}), we assign the velocity at time step $n$ as
\begin{equation}
	\mathbf{U}_i^n = -\frac{\partial E}{\partial \mathbf{x}_i},
\end{equation}
where the right-hand side can be calculated from Eq. (\ref{eq:ngws}) for interior domain and Eq. (\ref{eq:ngb}) for bounded domains.
With the velocity $\mathbf{U}_i^n$, we can update the particle position by Eq. (\ref{eq:stokes}), with $c=0.01$.

We would like to stress that, because the constant coefficient in $\mathbf{U}_i^n$ is offset by $\max |\mathbf{U}_i^n|$ in Eq. (\ref{eq:controlsize}), we can replace $\tilde{v}$
in Eqs. (\ref{eq:ngws}) and (\ref{eq:ngb}) by $v_0 \approx (\Delta x)^d$, without loss of generality. 
\section{Numerical validation of relaxed particle distribution patterns}
\label{validationpatt}
If we let the real volume of all particles be unity, there is a unique distribution for square and hexagonal lattices respectively. However, there can be infinitely many possible distributions for the other three types.
We introduce some parameters to represent different distributions for these three types.

\begin{figure}[tb!]
	\centering
	\includegraphics[width = 0.7\textwidth]{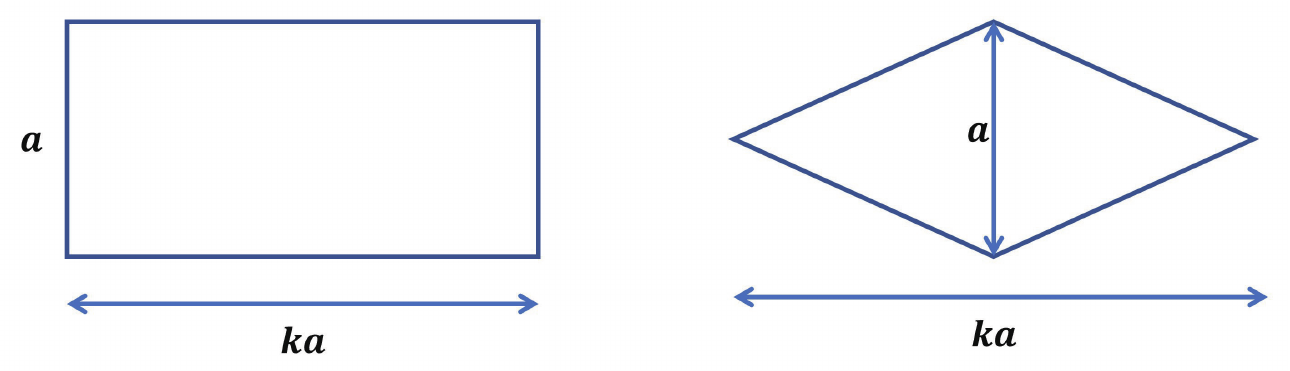}
	\caption{Determination of rectangular lattice (left panel) and diamond lattice (right panel) by the ratio $k$.}
	\label{fig:diamond}
\end{figure}

\begin{figure}[tb!]
	\centering
	\includegraphics[width = 0.3\textwidth]{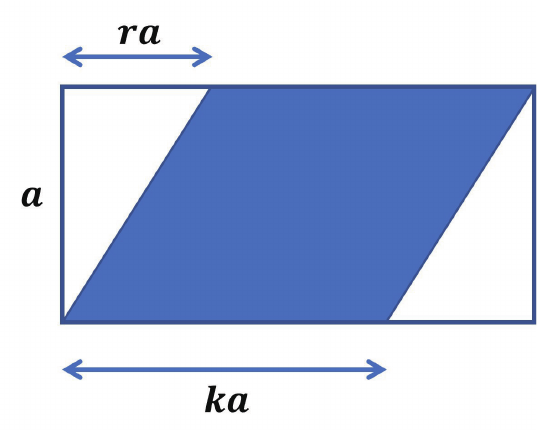}
	\caption{Determination of parallelogram lattice by the ratio $k$ and $r$.
	}
	\label{fig:parall}
\end{figure}

For rectangular and diamond distributions with unit volume, we only need two parameters, aspect ratio $k$ and the length of the shorter edge/axis, to determine its shape, see Fig. \ref{fig:diamond}.
For parallelogram distributions, we need three parameters including the height $a$, the ratio $k$ of base and height, and the 
degree of slant $r$, see Fig. \ref{fig:parall}.
We can find the largest characteristic volume of each type.
For example, the largest characteristic volume of rectangular distributions is obtained by
\begin{equation}
v_{rectangle} = \max_{k,a} \left\{\frac{1}{\sum_{i,j} W(\sqrt{(i^2k^2+j^2)}a,h)} \right\}.
\end{equation}

\begin{figure}[htbp]
	\centering
	\subfigure[$h=1.5\Delta x$]{\includegraphics[width=0.45\textwidth]{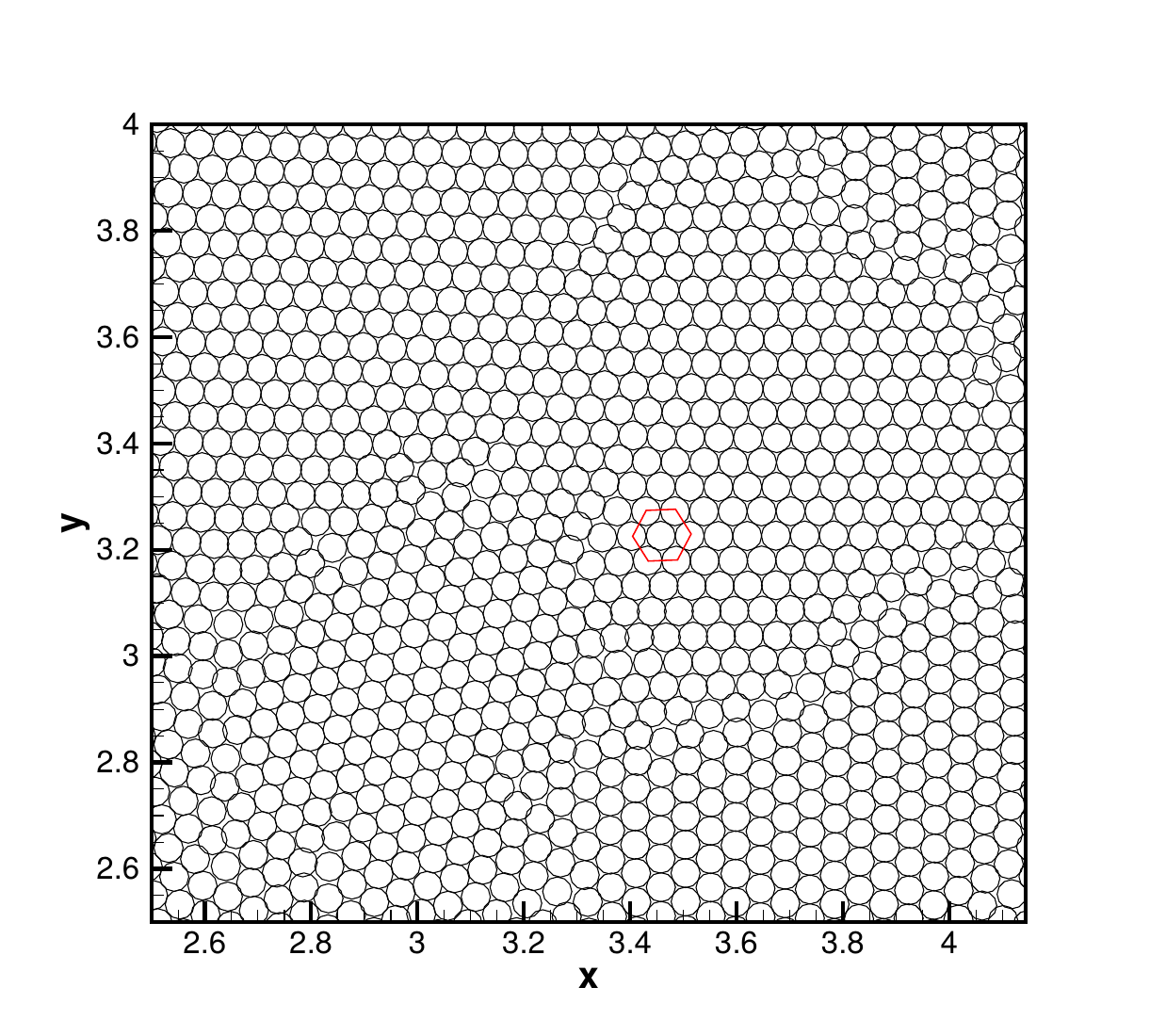}\label{fig:image1.5}}
	\hfill
	\subfigure[$h=1.8\Delta x$]{\includegraphics[width=0.45\textwidth]{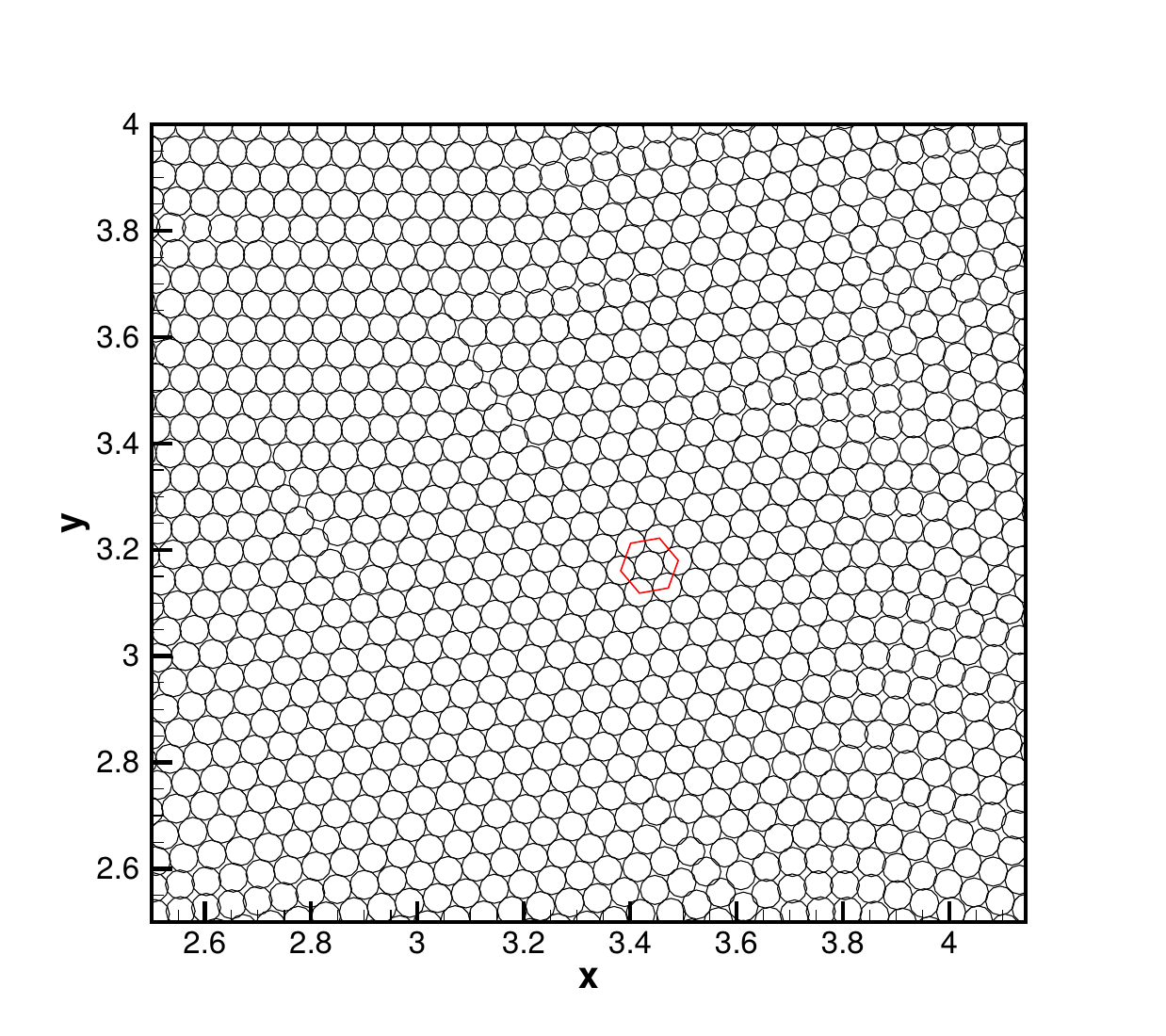}\label{fig:image1.8}}
	
	\vspace{0.5cm}
	\subfigure[$h=2.25\Delta x$]{\includegraphics[width=0.45\textwidth]{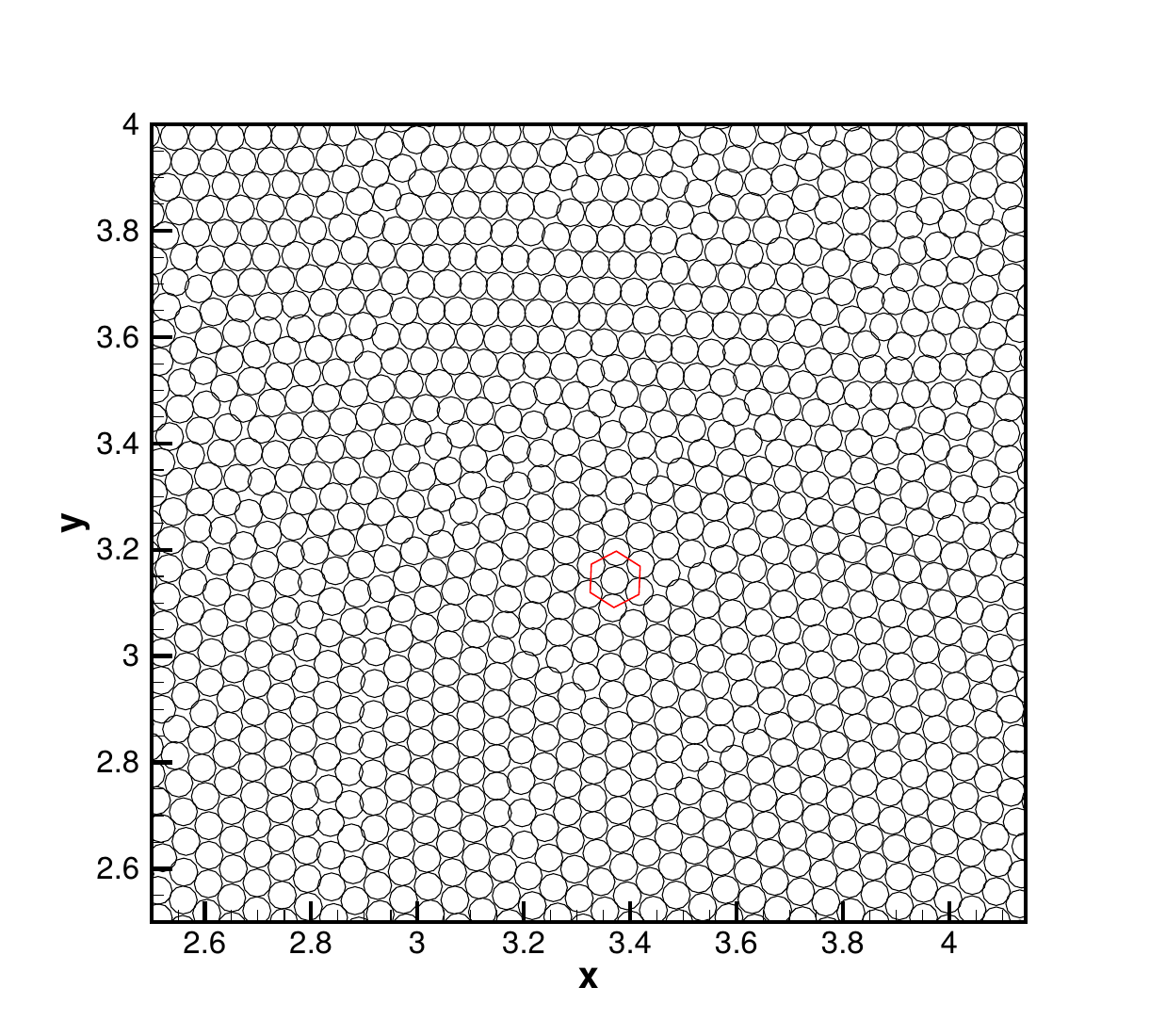}\label{fig:image2.25}}
	\hfill
	\subfigure[$h=2.75\Delta x$]{\includegraphics[width=0.45\textwidth]{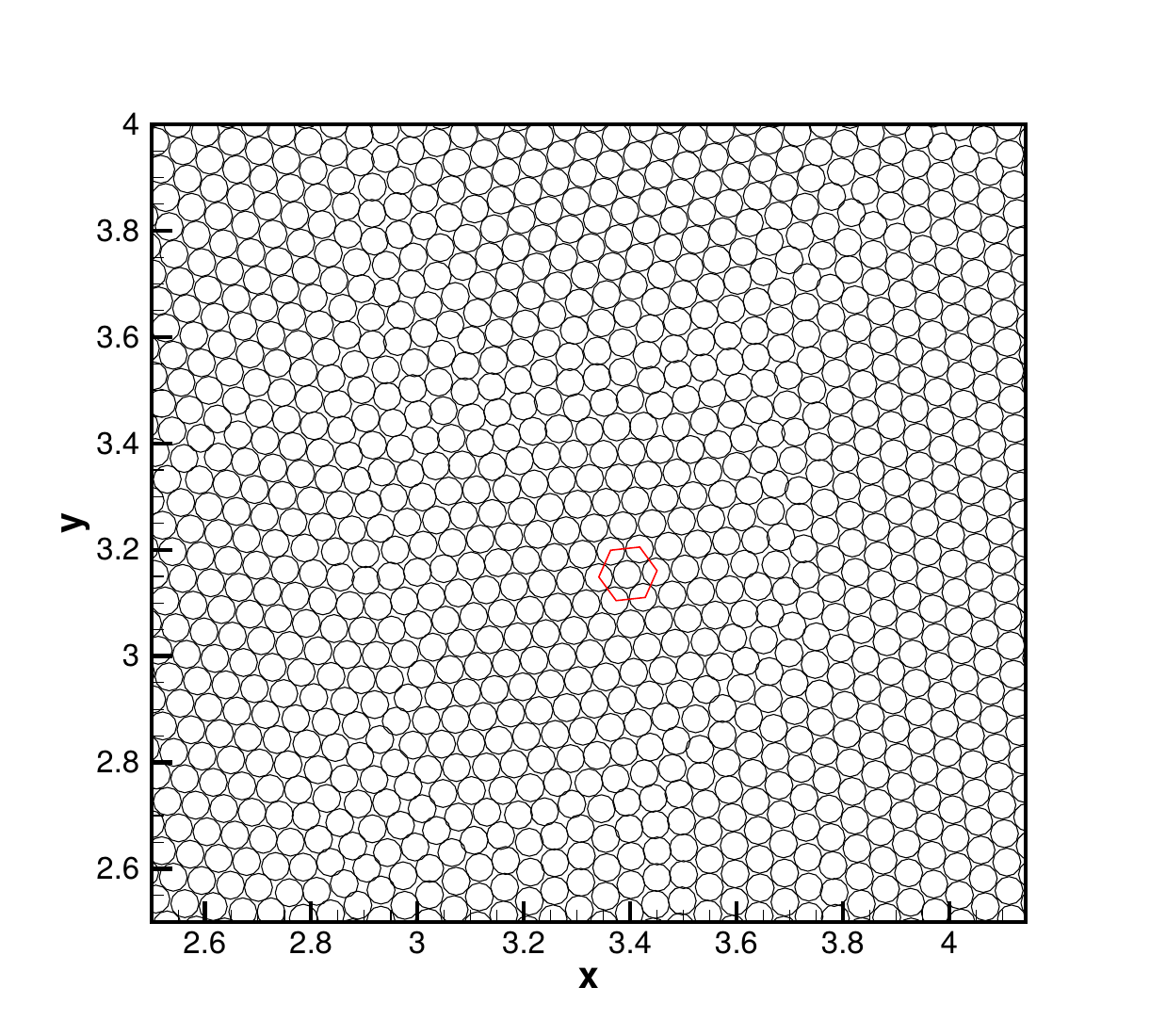}\label{fig:image2.75}}
	
	\vspace{0.5cm}
	\subfigure[$h=3.3\Delta x$]{\includegraphics[width=0.45\textwidth]{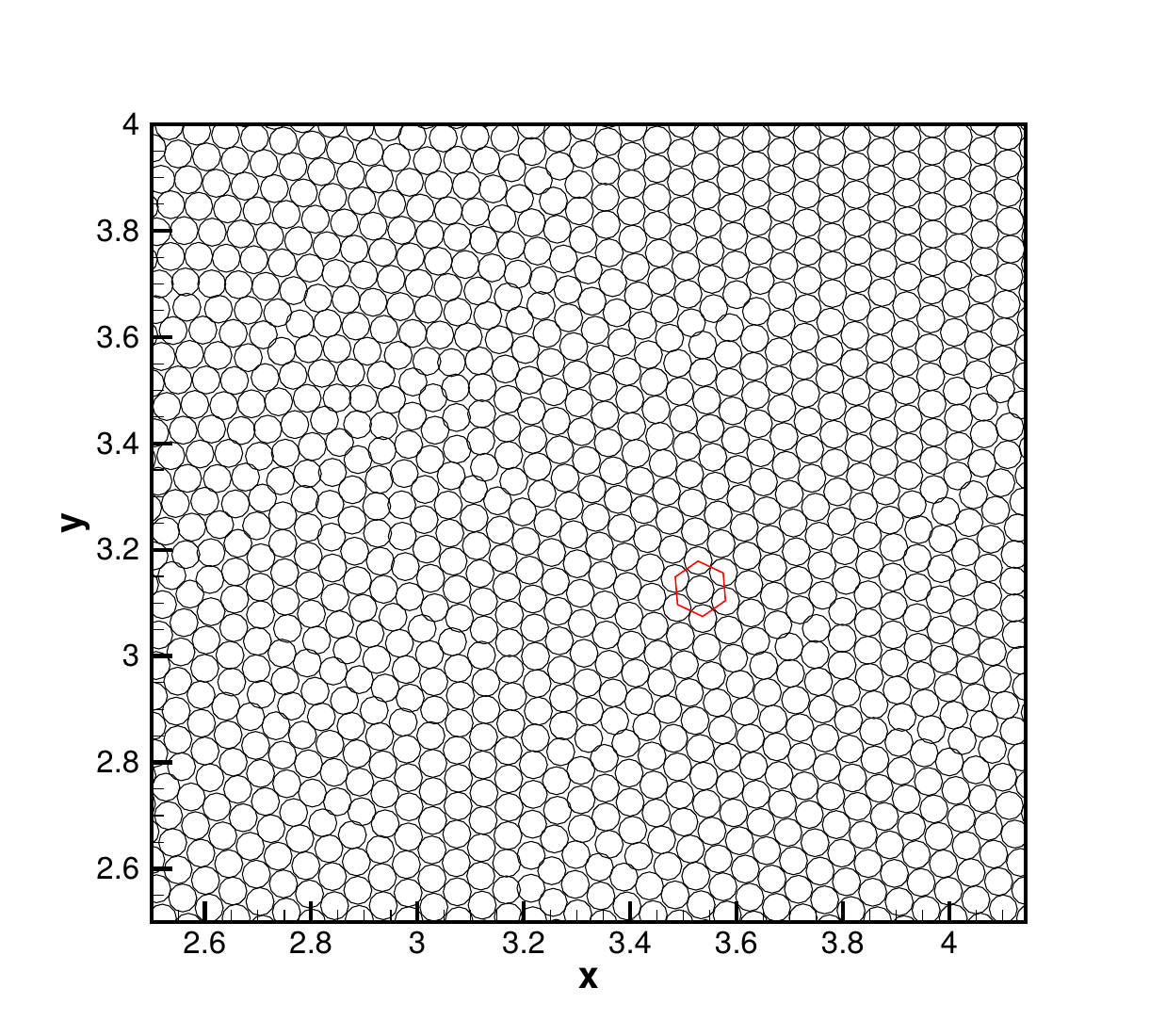}\label{fig:image3.3}}
	\hfill
	\subfigure[$h=3.45\Delta x$]{\includegraphics[width=0.45\textwidth]{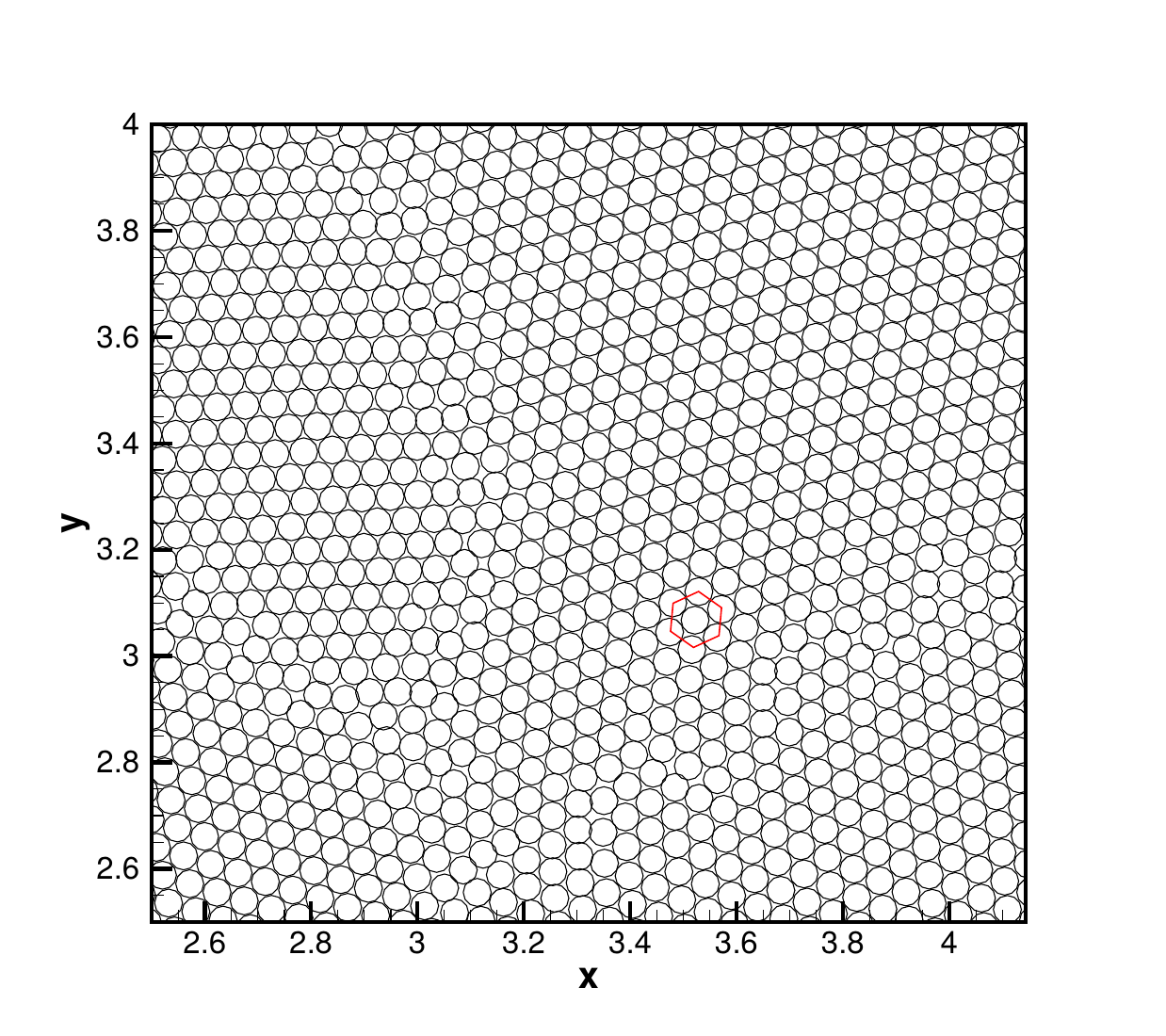}\label{fig:image3.45}}
	
	\caption{Hexagonal distributions in interior domain}
	\label{fig:hex}
\end{figure}

\begin{figure}[htbp]
	\centering
	\subfigure[$h=2.1\Delta x$]{\includegraphics[width=0.45\textwidth]{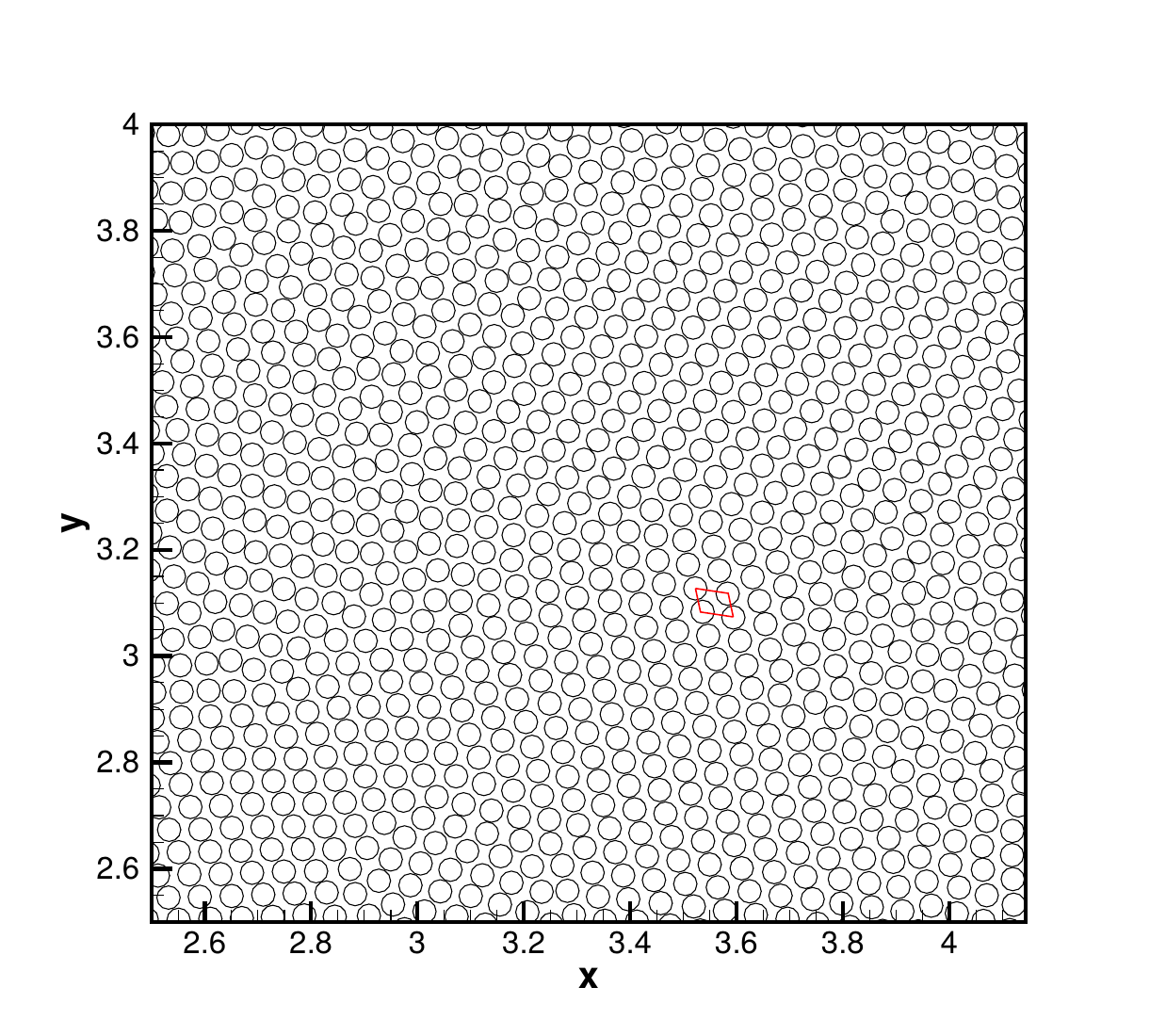}\label{fig:image2.1}}
	\hfill
	\subfigure[$h=2.2\Delta x$]{\includegraphics[width=0.45\textwidth]{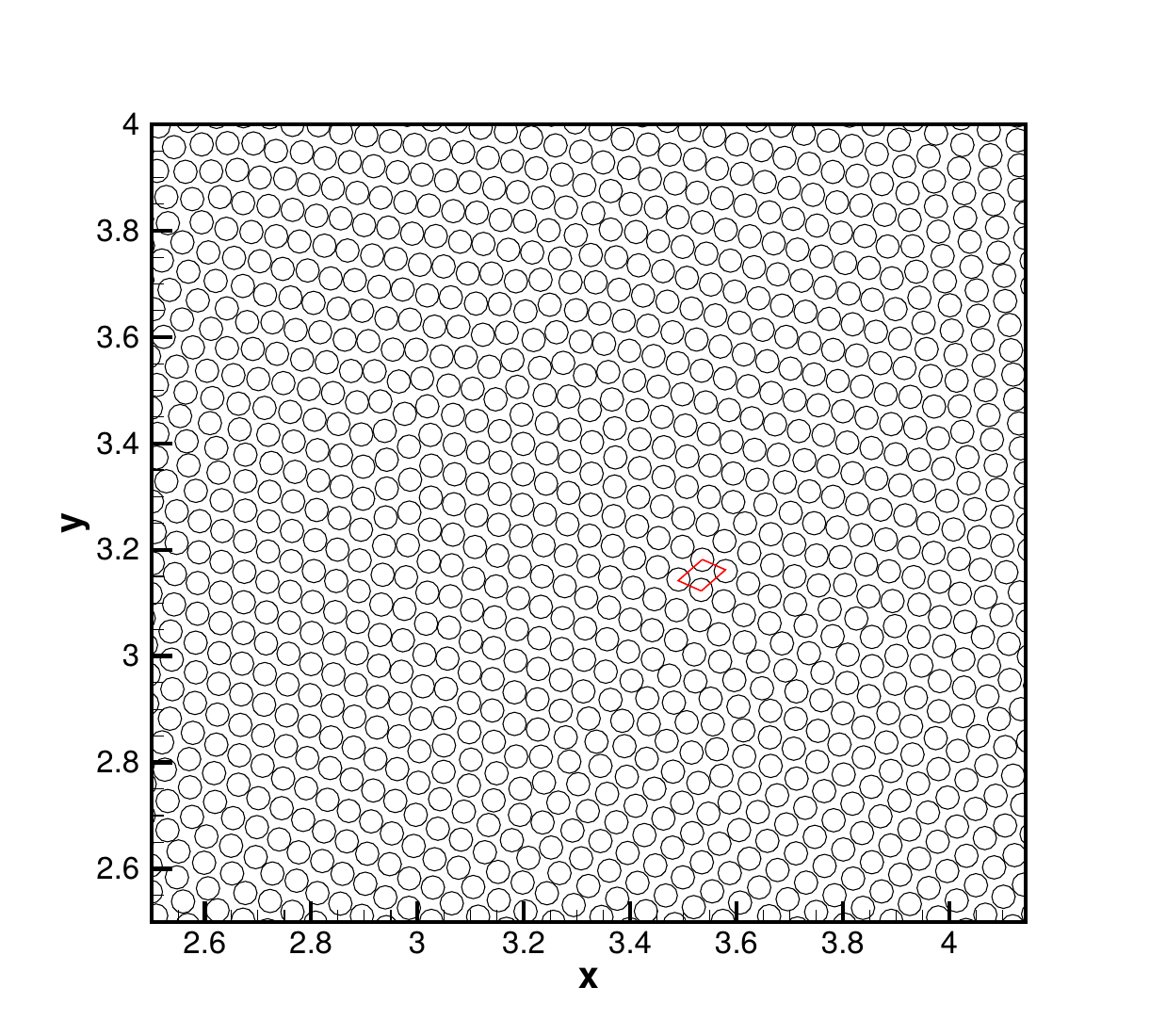}\label{fig:image2.2}}
	\vspace{0.5cm}
	\subfigure[$h=2.8\Delta x$]{\includegraphics[width=0.45\textwidth]{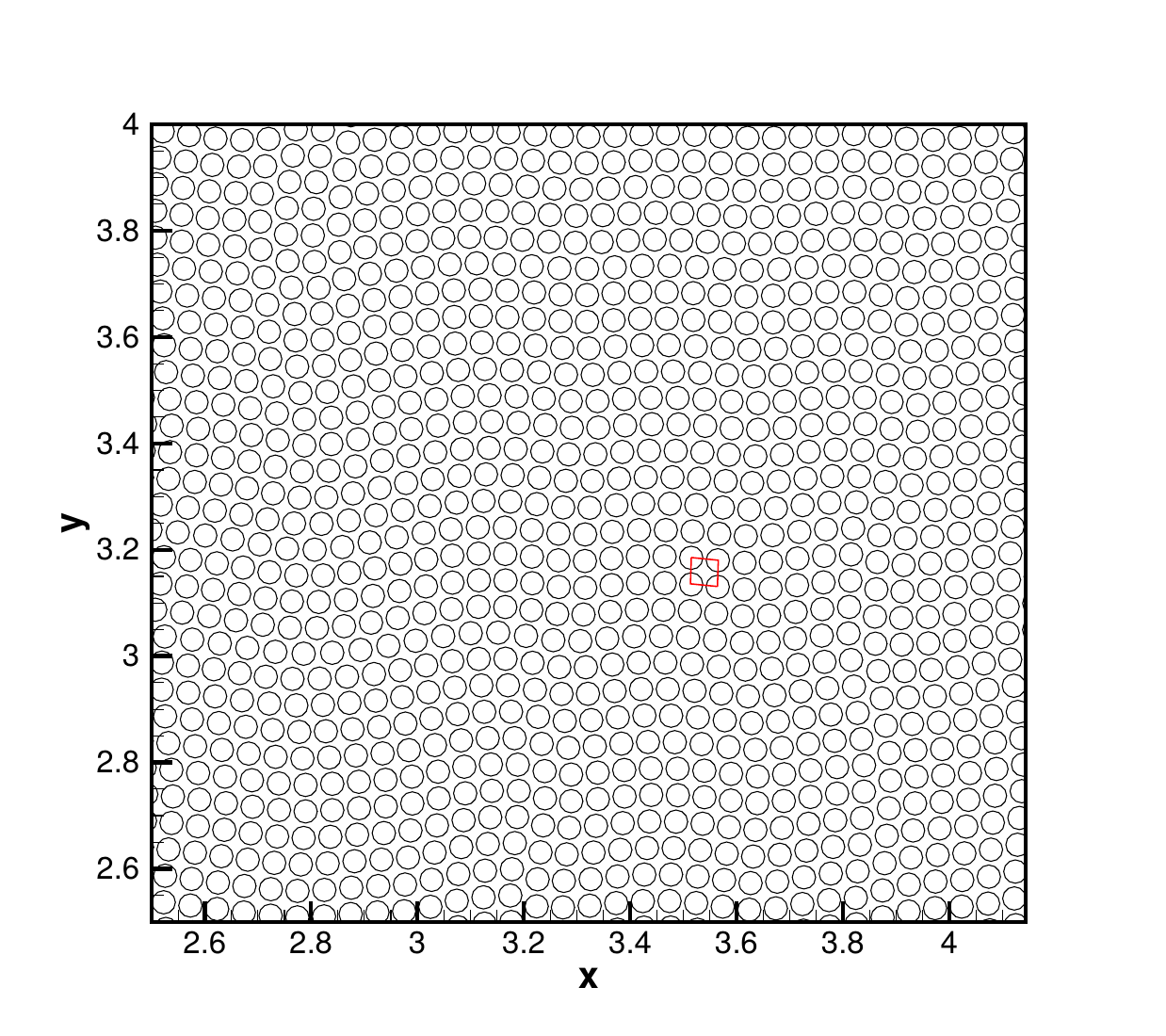}\label{fig:image2.8}}
	\hfill
	\subfigure[$h=3.25\Delta x$]{\includegraphics[width=0.45\textwidth]{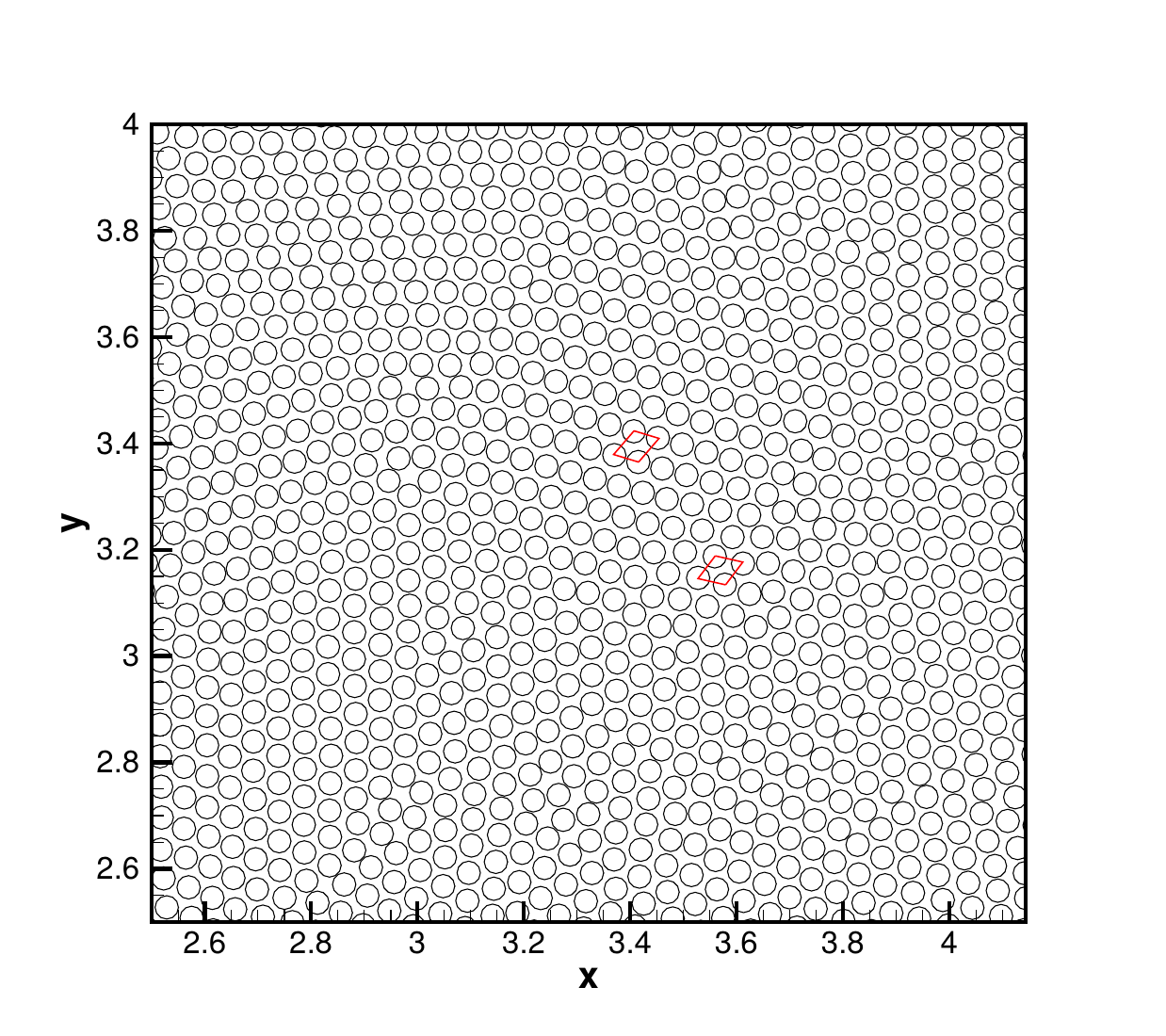}\label{fig:image3.25}}
	\caption{Parallelogram distributions in interior domain}
	\label{fig:para}
\end{figure}

\begin{figure}[htbp]
	\vspace{0.5cm}
	\subfigure[$h=1.85\Delta x$]{\includegraphics[width=0.45\textwidth]{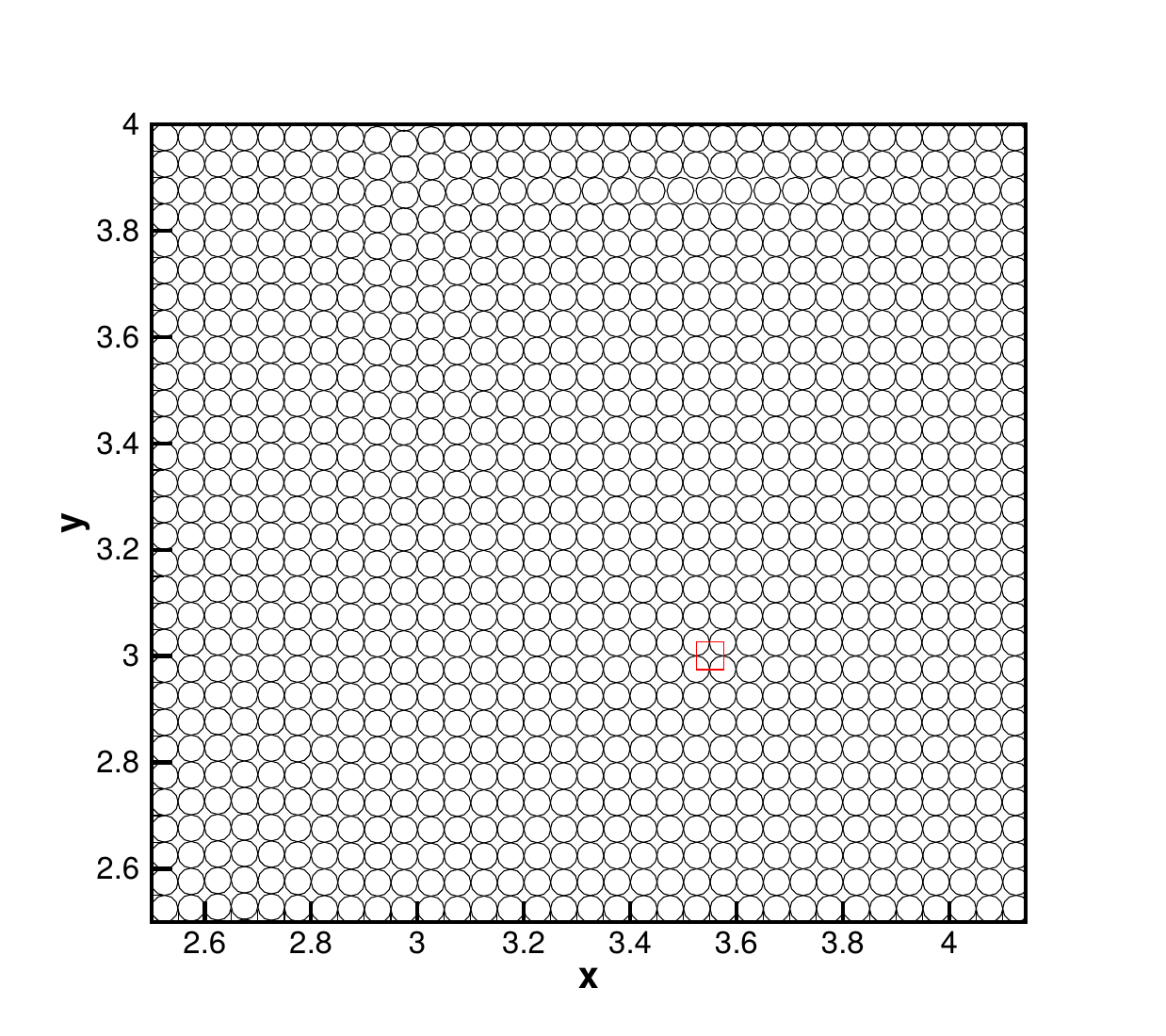}\label{fig:image1.85}}
	\hfill
	\subfigure[$h=2.05\Delta x$]{\includegraphics[width=0.45\textwidth]{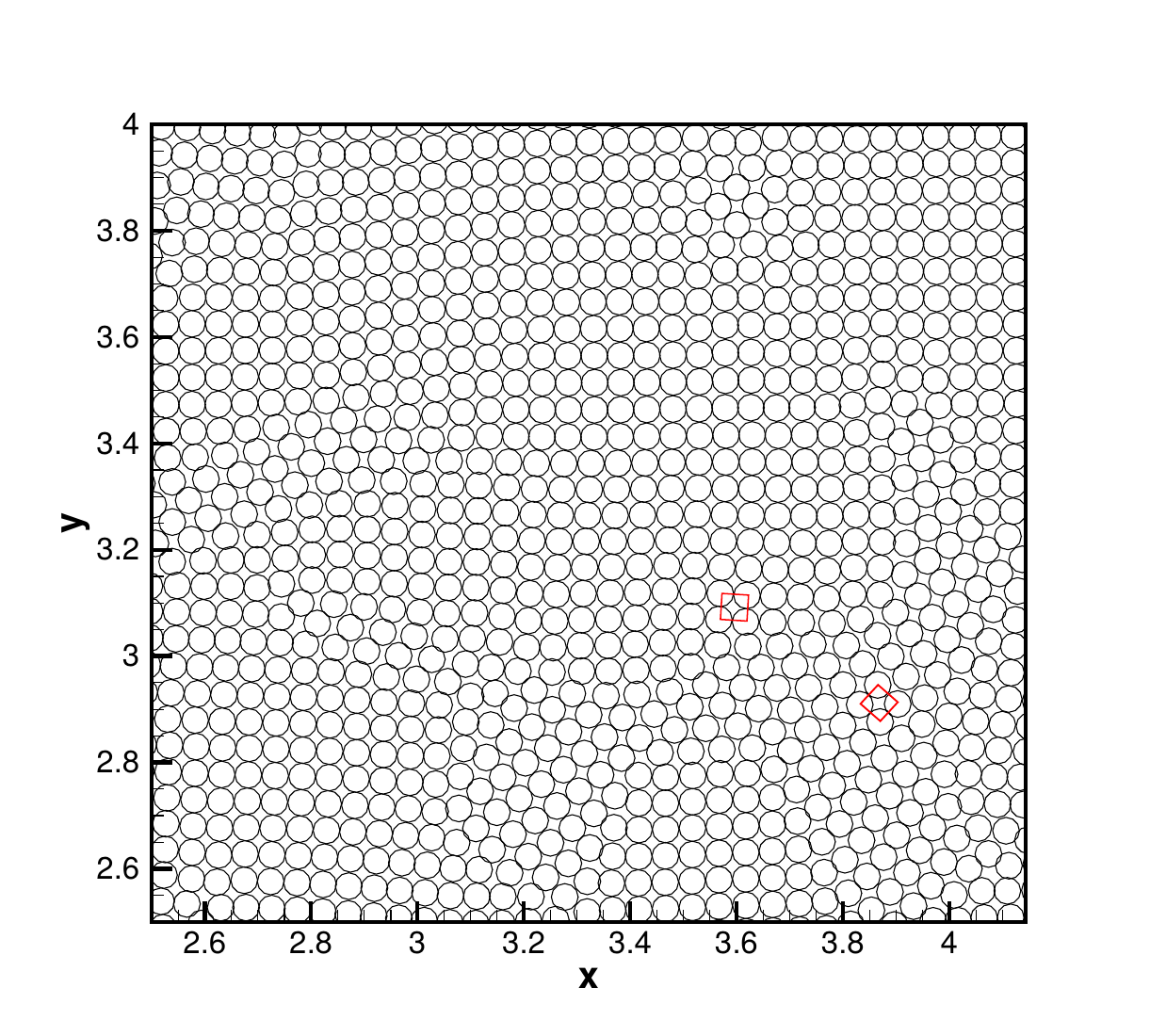}\label{fig:image2.05}}
	
	\caption{Square distributions in interior domain}
	\label{fig:square}
\end{figure}

\begin{figure}[htbp]
		\subfigure[$h=2.25dp$]{\includegraphics[width=0.45\textwidth]{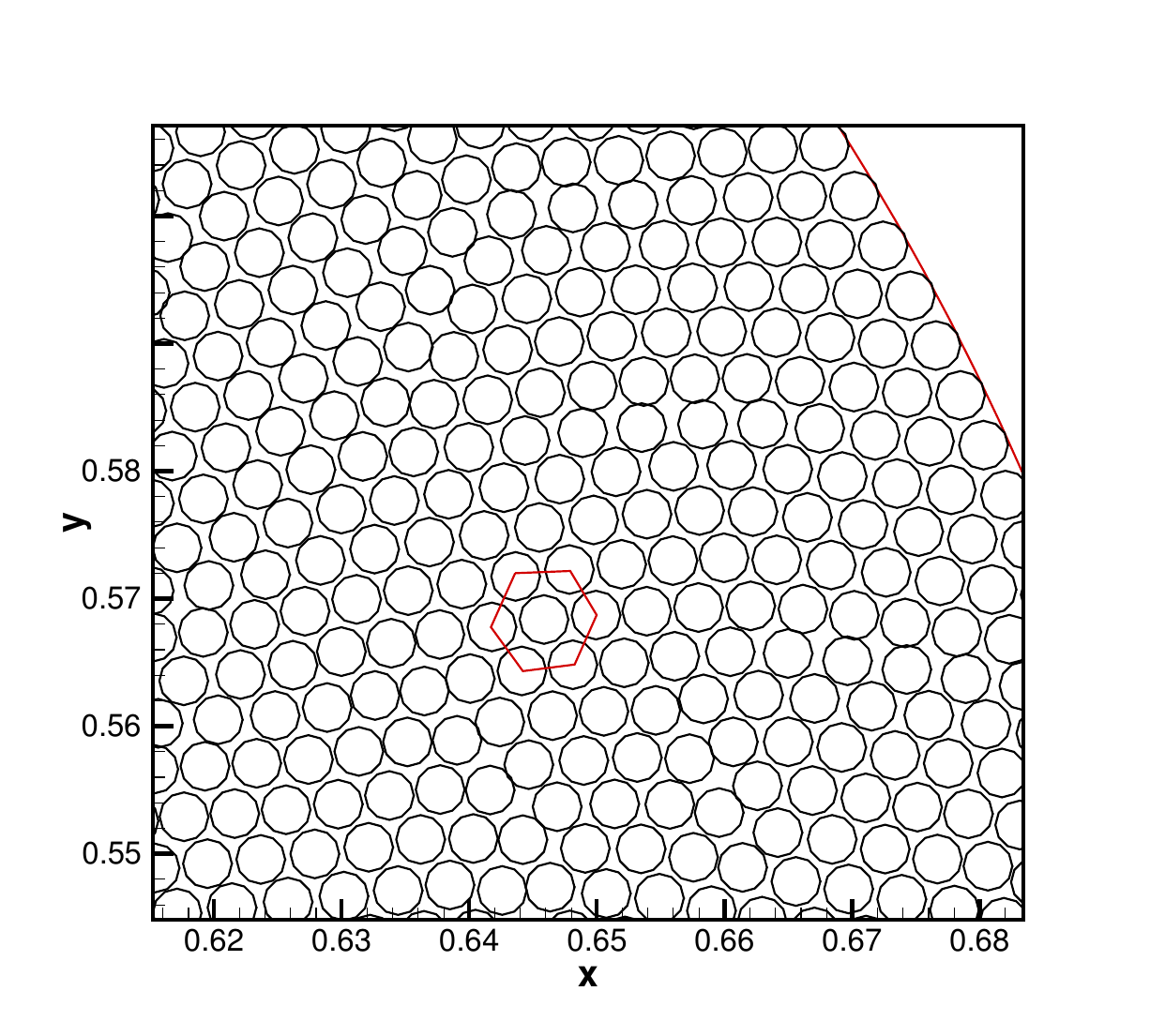}\label{fig:image225b}}
		\hfill
		\subfigure[$h=2.75dp$]{\includegraphics[width=0.45\textwidth]{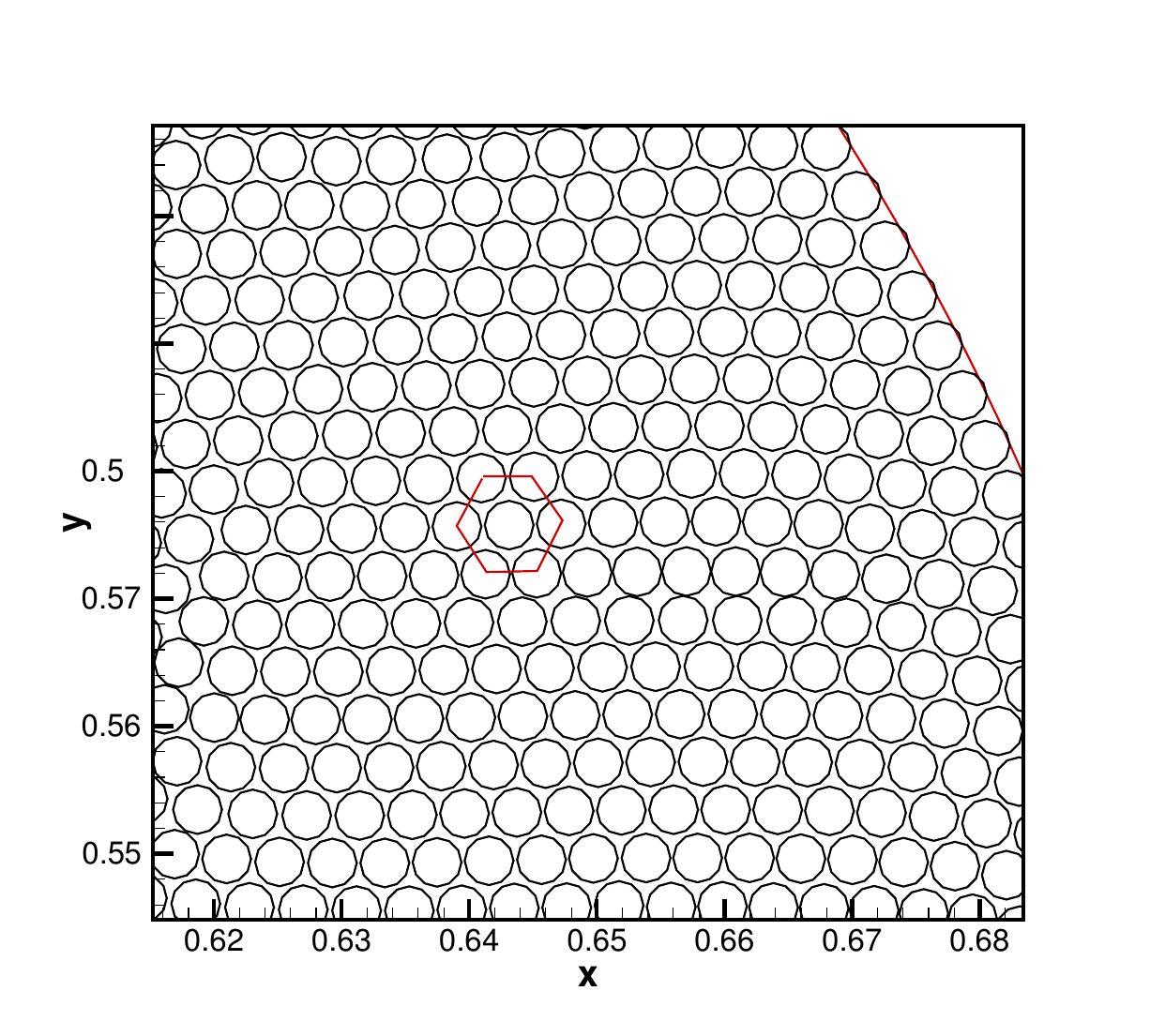}\label{fig:image275b}}
		\vspace{0.5cm}
		\subfigure[$h=3.3dp$]{\includegraphics[width=0.45\textwidth]{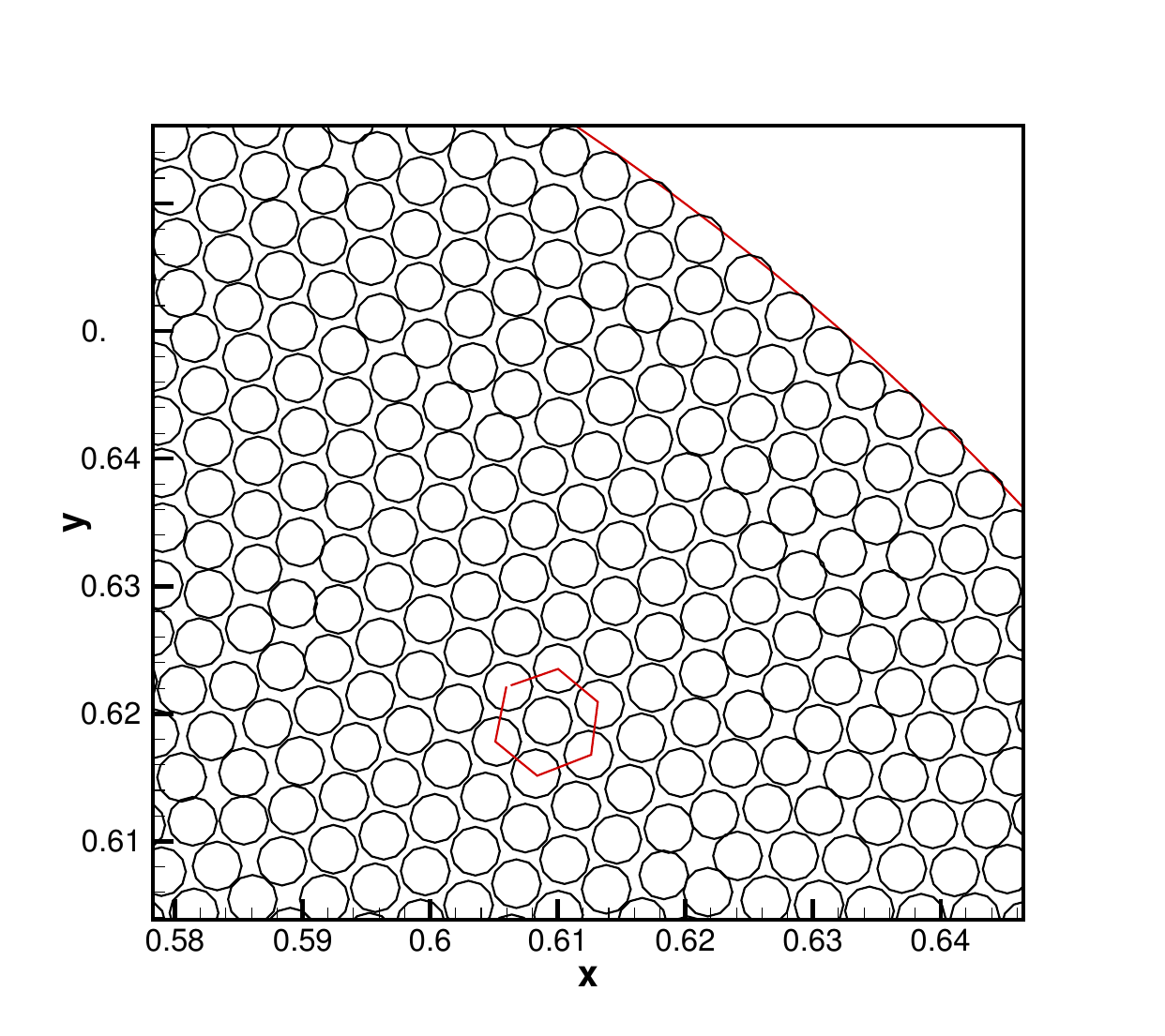}\label{fig:image330b}}
		\hfill
		\subfigure[$h=3.5dp$]{\includegraphics[width=0.45\textwidth]{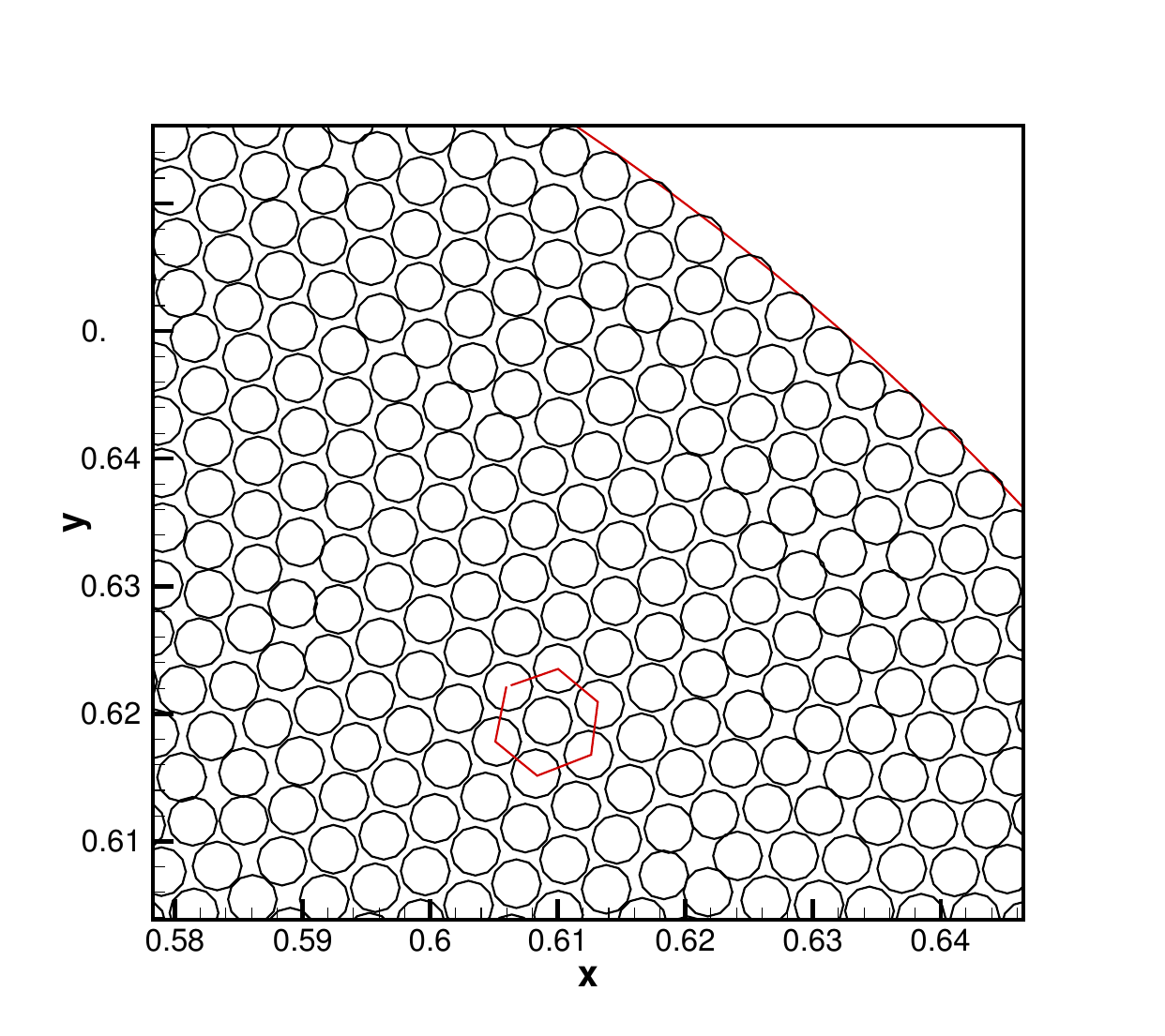}\label{fig:image350b}}
		
		\caption{Hexagonal distributions in a circular bounded domain}
		\label{fig:hexb}
	\end{figure}

\begin{figure}[htbp]
		\subfigure[$h=2.0dp$]{\includegraphics[width=0.45\textwidth]{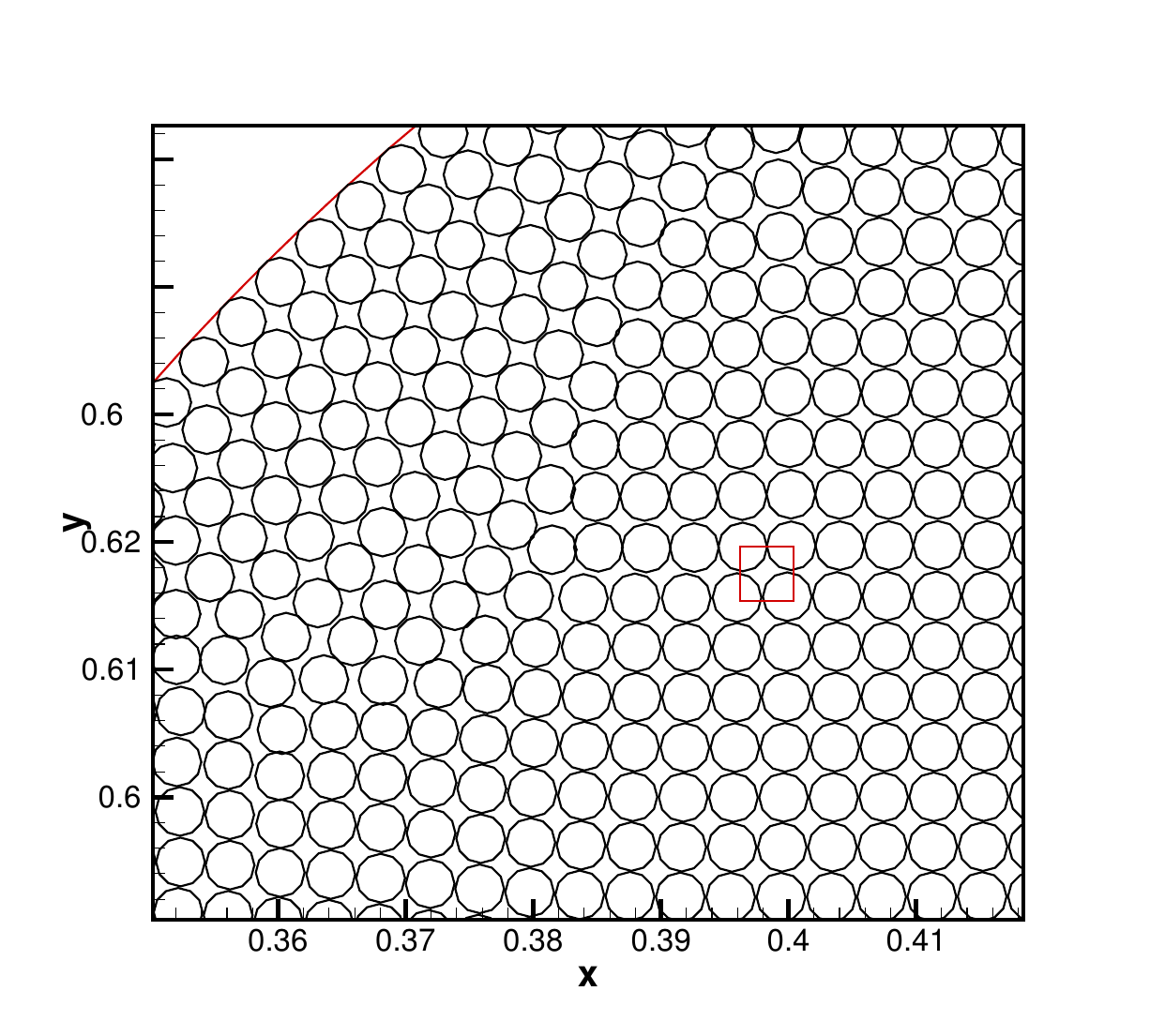}\label{fig:image200b}}
		\hfill
		\subfigure[$h=2.05dp$]{\includegraphics[width=0.45\textwidth]{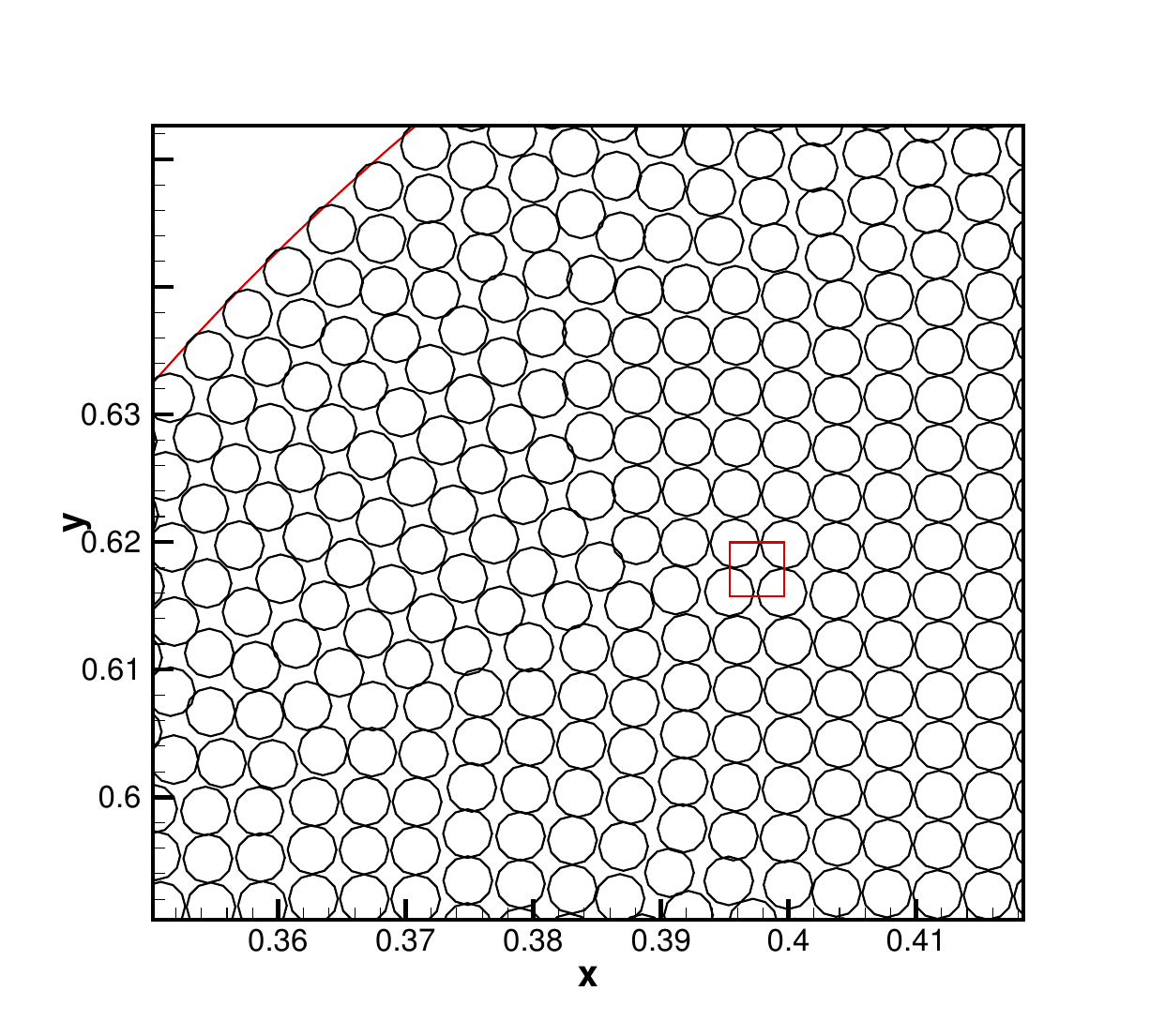}\label{fig:image205b}}
		\caption{Square distributions in a circular bounded domain}
		\label{fig:squareb}
	\end{figure}

\begin{figure}[htbp]
		\subfigure[$h=2.1dp$]{\includegraphics[width=0.45\textwidth]{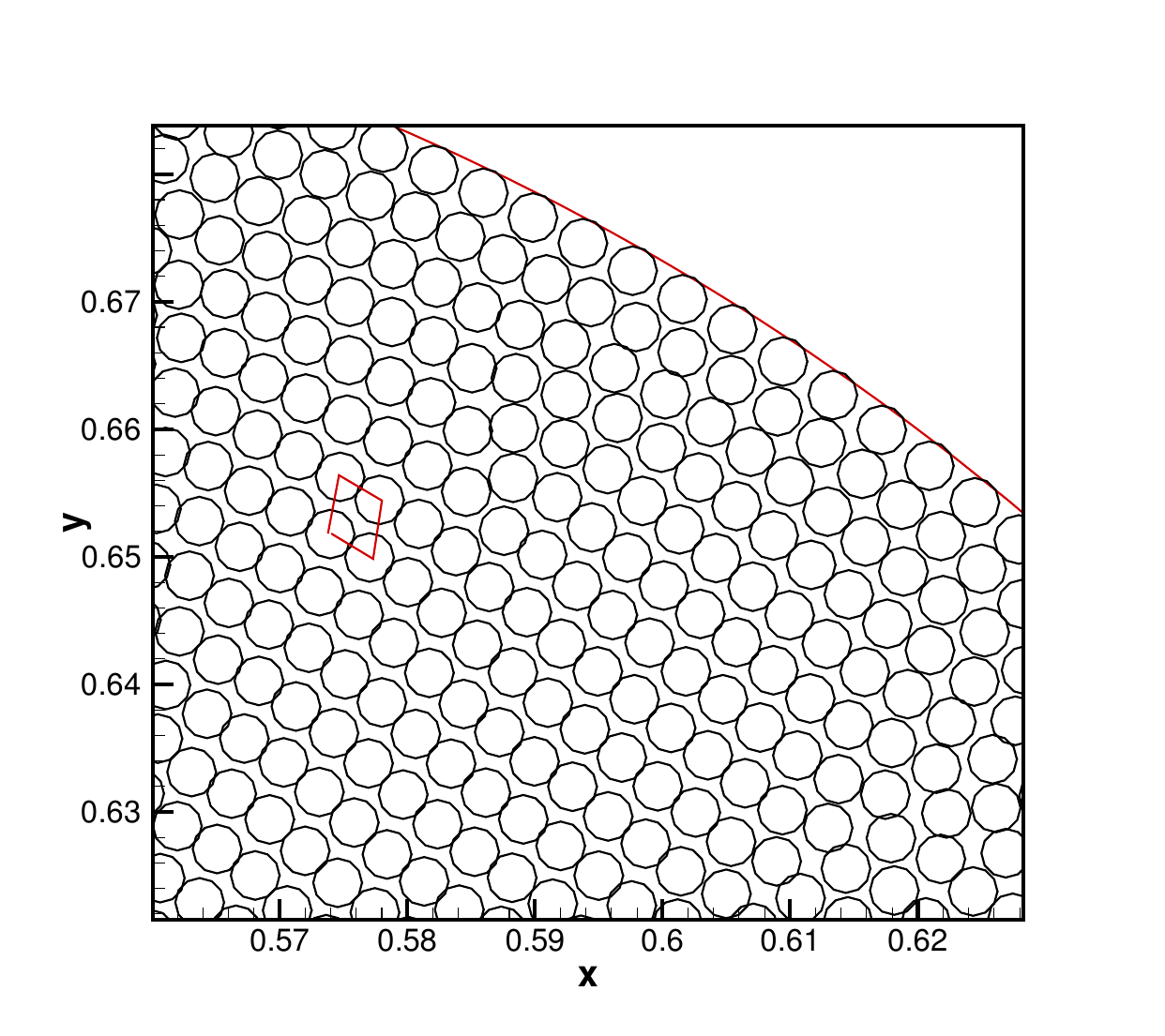}\label{fig:image210b}}
		\hfill
		\subfigure[$h=2.2dp$]{\includegraphics[width=0.45\textwidth]{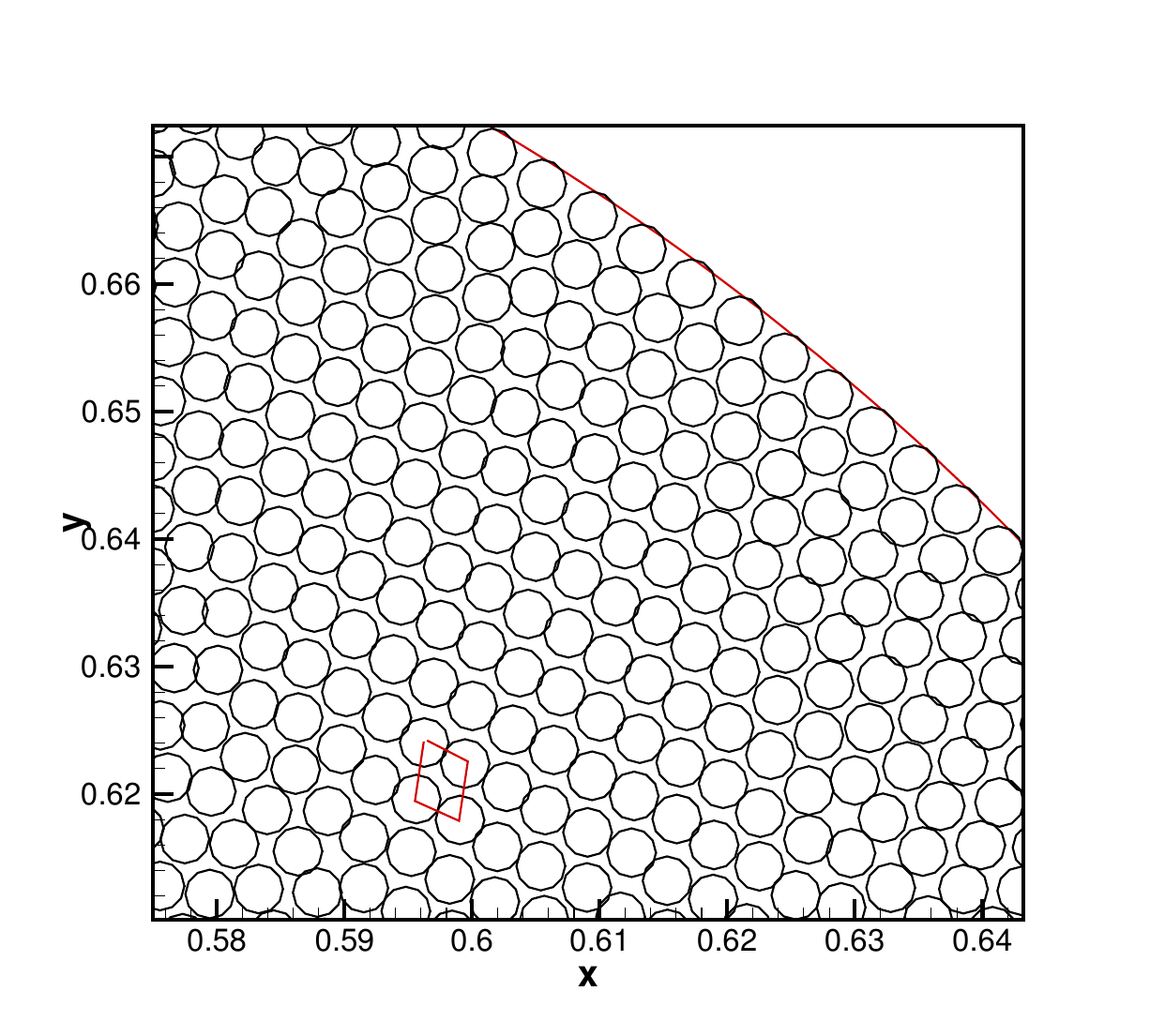}\label{fig:image220b}}
		\vspace{0.5cm}
		\subfigure[$h=2.8dp$]{\includegraphics[width=0.45\textwidth]{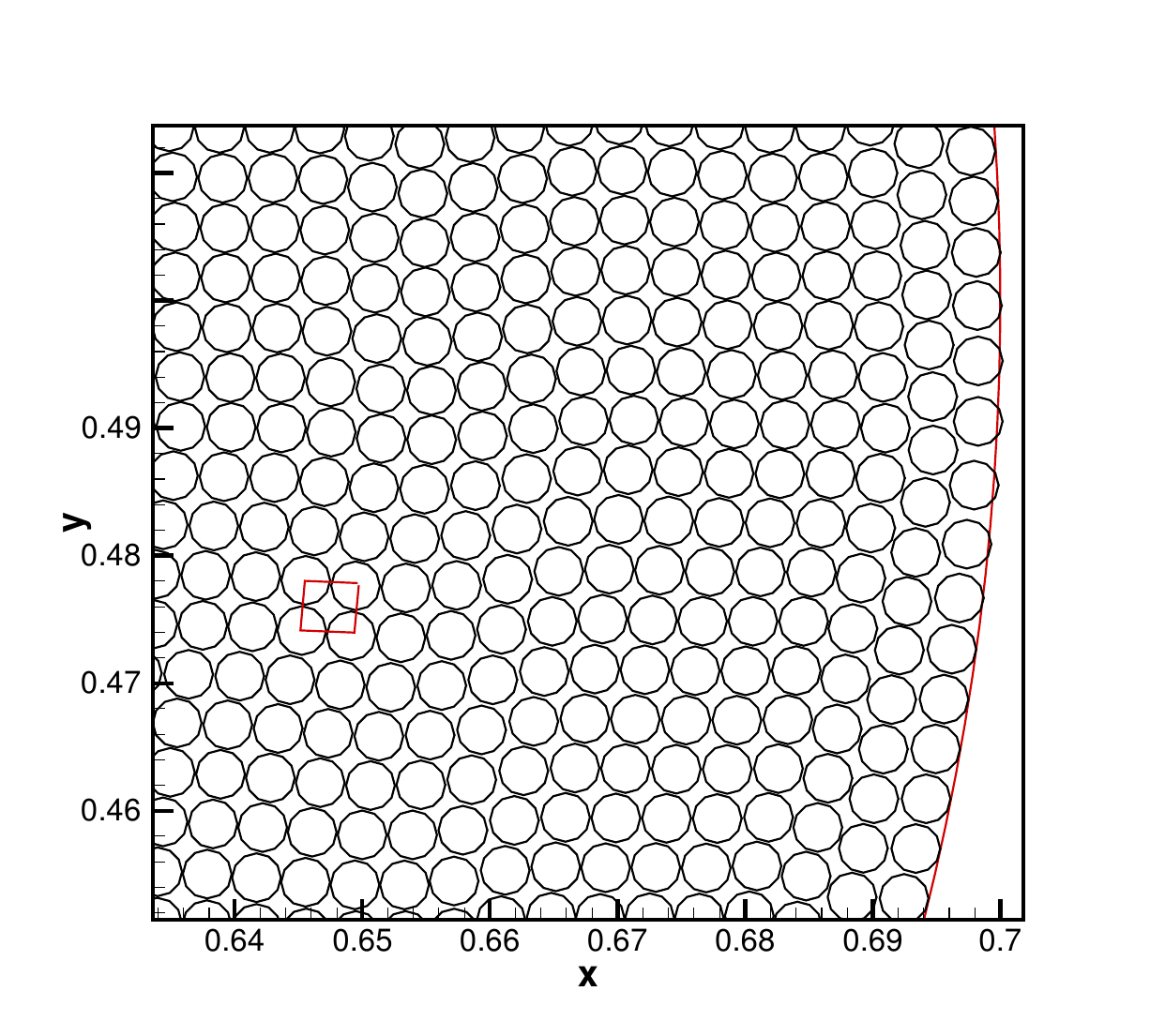}\label{fig:image280b}}
		%\label{combination}
		\hfill
		\subfigure[$h=3.25dp$]{\includegraphics[width=0.45\textwidth]{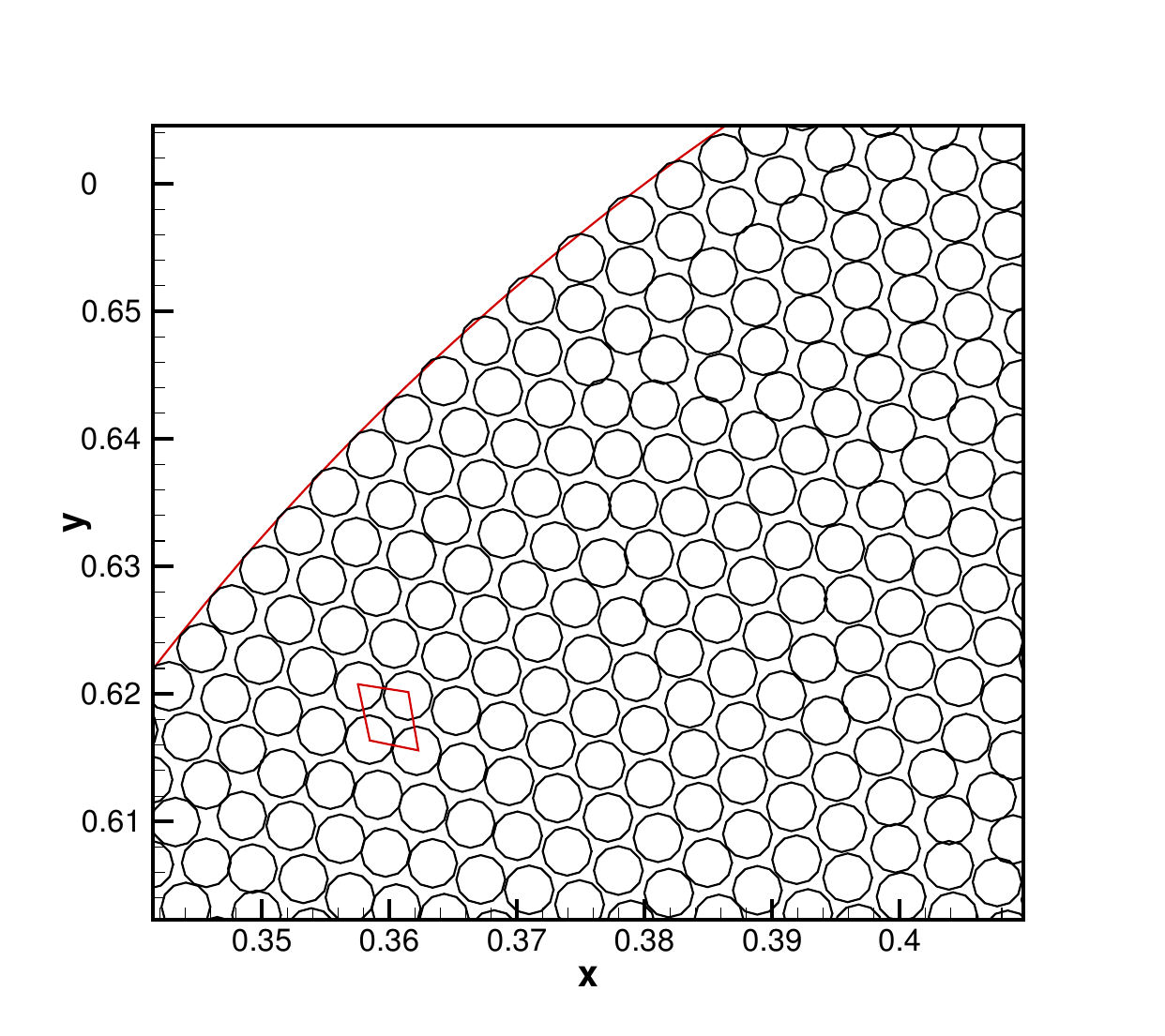}\label{fig:image325b}}
		\caption{Parallelogram distributions in a circular bounded domain.}
		\label{fig:parab}
	\end{figure}

Tab. \ref{tab:fullprediction} shows the largest characteristic volume of different types and the predicted patterns, which corresponds to the largest characteristic volume of all types.
Fig. \ref{fig:hex}, Fig. \ref{fig:para}, and Fig. \ref{fig:square} show the distribution near the corresponding critical values of hexagonal, parallelogram, and square distributions respectively.
We also verify our prediction in bounded domains such as a circle.
The corresponding patterns are shown in Fig. \ref{fig:hexb}, Fig. \ref{fig:squareb}, and Fig. \ref{fig:parab}.
%We need to stress that the particle relaxation method not always produces perfect patterns.
%First, the size of the domain may not be compatible with the predicted pattern.
%Take the square lattice as an example, whose $|\mathbf{e}_1|=|\mathbf{e}_2|=1$, but the size of the domain can have a different ratio.
%Secondly, the accuracy of the particle relaxation method may not be sufficient for the predicted perfect pattern because the difference in the energy between the predicted pattern and other patterns is small.
%This can be improved by using a smaller time-step size. 
\begin{table}[tb!]
	\centering
	\caption{Details of Predictions for 2D Lattice}
	%\resizebox{.5\textwidth}{\textwidth}{!}{
		\begin{adjustbox}{width=\textwidth,center}
			\begin{tabular}{|c|c|c|c|c|c|c|}
				\hline
				$h/\Delta x$ & $v_{\text{hexagon}}$ & $v_{\text{square}}$ & $v_{\text{diamond}}$ & $v_{\text{rectangle}}$ & $v_{\text{parallelogram}}$ & Pattern(largest characteristic volume) \\ \hline
				1.5 & 0.878 & 0.8548 & 0.872 & 0.8548 & 0.872 & Hexagon \\ \hline
				1.55 & 0.8982 & 0.8779 & 0.8927 & 0.8779 & 0.8927 & Hexagon \\ \hline
				1.6 & 0.9138 & 0.8975 & 0.909 & 0.8975 & 0.909 & Hexagon \\ \hline
				1.65 & 0.9257 & 0.9137 & 0.9217 & 0.9137 & 0.9217 & Hexagon \\ \hline
				1.7 & 0.9346 & 0.9269 & 0.9314 & 0.9269 & 0.9314 & Hexagon \\ \hline
				1.75 & 0.9413 & 0.9373 & 0.9389 & 0.9373 & 0.9389 & Hexagon \\ \hline
				1.8 & 0.9465 & 0.9455 & 0.9449 & 0.9455 & 0.9455 & Hexagon \\ \hline
				1.85 & 0.9508 & 0.9519 & 0.9499 & 0.9519 & 0.9519 & Square \\ \hline
				1.9 & 0.9546 & 0.9568 & 0.9543 & 0.9568 & 0.9568 & Square \\ \hline
				1.95 & 0.9582 & 0.9606 & 0.9583 & 0.9606 & 0.9606 & Square \\ \hline
				2.0 & 0.9616 & 0.9638 & 0.9621 & 0.9638 & 0.9638 & Square \\ \hline
				2.05 & 0.965 & 0.9664 & 0.966 & 0.9664 & 0.9664 & Square \\ \hline
				2.1 & 0.9684 & 0.9688 & 0.9695 & 0.9689 & 0.9696 & Parallelogram, k = 1.18, r = 0.42 \\ \hline
				2.15 & 0.9716 & 0.9711 & 0.9726 & 0.9714 & 0.9727 & Parallelogram, k = 1.19, r = 0.42 \\ \hline
				2.2 & 0.9747 & 0.9734 & 0.9752 & 0.9737 & 0.9753 & Parallelogram, k = 1.17, r = 0.42 \\ \hline
				2.25 & 0.9777 & 0.9756 & 0.9777 & 0.9758 & 0.9777 & Hexagon \\ \hline
				2.25 & 0.9777 & 0.9756 & 0.9777 & 0.9758 & 0.9777 & Hexagon \\ \hline
				2.3 & 0.9805 & 0.9779 & 0.9803 & 0.9779 & 0.9803 & Hexagon \\ \hline
				2.35 & 0.9831 & 0.9801 & 0.9826 & 0.9801 & 0.9826 & Hexagon \\ \hline
				2.4 & 0.9854 & 0.9824 & 0.9847 & 0.9824 & 0.9847 & Hexagon \\ \hline
				2.45 & 0.9873 & 0.9845 & 0.9865 & 0.9845 & 0.9865 & Hexagon \\ \hline
				2.5 & 0.9889 & 0.9864 & 0.9881 & 0.9864 & 0.9881 & Hexagon \\ \hline
				2.55 & 0.9902 & 0.9881 & 0.9894 & 0.9881 & 0.9894 & Hexagon \\ \hline
				2.6 & 0.9912 & 0.9896 & 0.9904 & 0.9896 & 0.9904 & Hexagon \\ \hline
				2.65 & 0.992 & 0.9909 & 0.9913 & 0.9909 & 0.9913 & Hexagon \\ \hline
				2.7 & 0.9926 & 0.9919 & 0.992 & 0.9919 & 0.992 & Hexagon \\ \hline
				2.75 & 0.993 & 0.9927 & 0.9926 & 0.9927 & 0.9927 & Hexagon \\ \hline
				2.8 & 0.9934 & 0.9934 & 0.9931 & 0.9934 & 0.9934 & Parallelogram, k = 1.0, r = 0.86 \\ \hline
				2.85 & 0.9937 & 0.9939 & 0.9936 & 0.9939 & 0.9939 & Parallelogram, k = 1.0, r = 0.82 \\ \hline
				2.9 & 0.994 & 0.9943 & 0.994 & 0.9943 & 0.9944 & Parallelogram, k = 1.0, r = 0.19 \\ \hline
				2.95 & 0.9943 & 0.9946 & 0.9944 & 0.9946 & 0.9947 & Parallelogram, k = 1.0, r = 0.79 \\ \hline
				3.0 & 0.9946 & 0.9949 & 0.9949 & 0.9949 & 0.995 & Parallelogram, k = 1.0, r = 0.23 \\ \hline
				3.05 & 0.995 & 0.9951 & 0.9953 & 0.9951 & 0.9953 & Parallelogram, $k=1.11$, $r=0.53$ \\ \hline
				3.1 & 0.9954 & 0.9954 & 0.9957 & 0.9954 & 0.9957 & Parallelogram, $k=1.14$, $r=0.44$ \\ \hline
				3.15 & 0.9958 & 0.9956 & 0.996 & 0.9957 & 0.9961 & Parallelogram, $k=1.15$, $r=0.44$ \\ \hline
				3.2 & 0.9962 & 0.9959 & 0.9963 & 0.9961 & 0.9963 & Parallelogram, $k=1.15$, $r=0.44$ \\ \hline
				3.25 & 0.9965 & 0.9961 & 0.9965 & 0.9963 & 0.9966 & Parallelogram, $k=1.12$, $r=0.44$ \\ \hline
				3.3 & 0.9969 & 0.9964 & 0.9968 & 0.9966 & 0.9968 & Hexagon\\ \hline
				3.35 & 0.9972 & 0.9967 & 0.9971 & 0.9968 & 0.9971 & Hexagon \\ \hline
				3.4 & 0.9975 & 0.997 & 0.9973 & 0.997 & 0.9973 & Hexagon \\ \hline
				3.45 & 0.9977 & 0.9972 & 0.9975 & 0.9972 & 0.9975 & Hexagon \\ \hline
				3.50 & 0.9979 & 0.9975 & 0.9977 & 0.9975 & 0.9977 & Hexagon \\ \hline
			\end{tabular}
		\end{adjustbox}
		\label{tab:fullprediction}
	\end{table}
	\section{Application to multi-region cases}
	\label{chap:nest}
\begin{figure}[H]
		\centering
		\includegraphics[width = 0.5\textwidth]{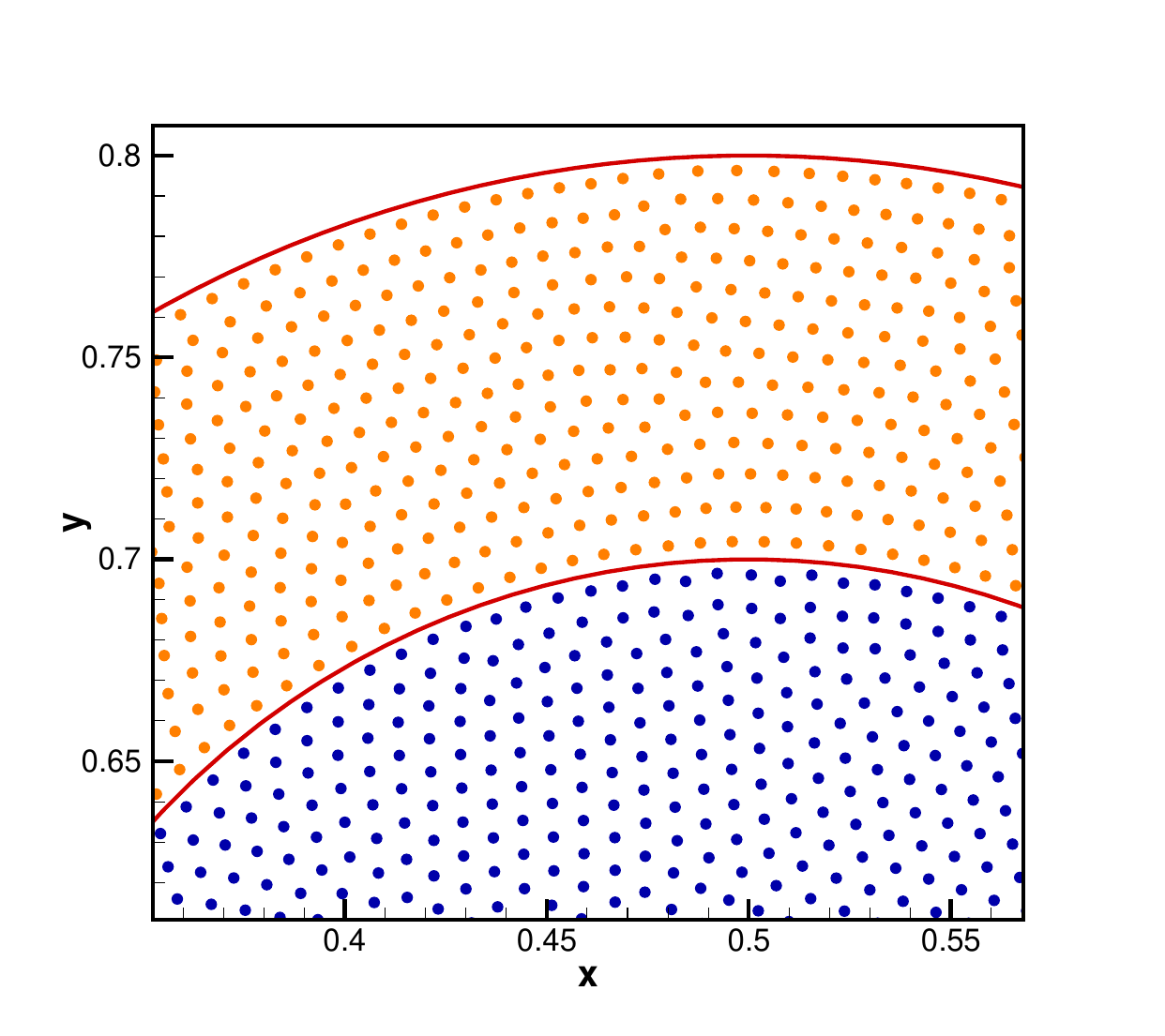}
		\caption{Two nested circles. The radii of them are $R_l=0.3$ and $R_s = 0.2$ respectively, with $128\times 128$ background mesh. The kernel cut-off radius is $h=2.6\Delta x$.
		}
		\label{fig:circles}
	\end{figure}

	We use a simple case to verify the formulation in Sec. \ref{chap:formula} for multiple bounded regions.
	In this case, we have two circles with radii $R_l=0.3$ and $R_s = 0.2$ respectively.
	We obtain the distributions in the inner circle and the annulus by the relaxation method independently, see Fig. \ref{fig:circles}.
	The particles in the inner circle and the annulus are distributed uniformly, and the distance between the boundary and boundary particles is well controlled.

\bibliographystyle{unsrt}
\bibliography{refs}
\end{document}

%% file: ex_shared.tex
\usepackage{lipsum}
\usepackage{amsfonts}
\usepackage{graphicx}
\usepackage{subfigure}
\usepackage{amsmath}
\usepackage{mathrsfs}
\usepackage{adjustbox}
\usepackage{algorithmic}
\usepackage[normalem]{ulem}
\ifpdf
  \DeclareGraphicsExtensions{.eps,.pdf,.png,.jpg}
\else
  \DeclareGraphicsExtensions{.eps}
\fi

\title{Analysis of the particle relaxation method for generating uniform particle distributions in smoothed particle hydrodynamics}

\author{Yu Fan\thanks{Corresponding Author, Chair of Aerodynamics and Fluid Mechanics, School of Engineering and Design, Technical University of Munich, Boltzmannstraße 15, 85748 Garching, Germany 
		(\email{yu.fan@tum.de})}, 
	\and Xiaoliang Li\thanks{Department of Mathematics ``Federigo Enriques", Universit\`a degli Studi di Milano, Via Cesare Saldini 50, 20133, Milan, Italy
	(\email{xiaoliang.li@unimi.it})}
	\and Shuoguo Zhang\thanks{Chair of Aerodynamics and Fluid Mechanics, School of Engineering and Design, Technical University of Munich, Boltzmannstraße 15, 85748 Garching, Germany
	(\email{shuoguo.zhang@tum.de}, \email{xiangyu.hu@tum.de}, \email{nikolaus.adams@tum.de}).}
	\and Xiangyu Hu\footnotemark[3]
	\and Nikolaus A. Adams\footnotemark[3]}
\usepackage{amsopn}